\documentclass[a4paper,twoside,12pt]{article}
\usepackage{graphicx,psfig,epsfig}
\usepackage{amsfonts}
\usepackage{amssymb}
\usepackage{amsmath}
\usepackage{subfigure}
\usepackage{setspace}
\usepackage[text={16cm,24cm},top=3cm,bottom=3cm,centering]{geometry}
 \usepackage{helvet,avant,bookman,mathptmx}
\renewcommand{\maketitle}%
{\thispagestyle{empty}
\noindent
\par\vspace*{2cm}\par
{\centering{\Large\textbf\MATCHtitle\par}
\par\vspace*{12pt}\par
{\textrm\MATCHauthor}\par\vspace*{12pt}\par}%
}
\renewenvironment{abstract}%
{\begin{center}\begin{minipage}[t]{1.0\textwidth}
\noindent\small\textbf{Abstract$\colon$}}%
{\end{minipage}\end{center}\par}

\newenvironment{keywords}%
{\begin{center}\begin{minipage}[t]{1.0\textwidth}
\noindent\small\textit{\bf Keywords:}}%
{\end{minipage}\end{center}\par}

\newenvironment{class}%
{\begin{center}\begin{minipage}[t]{1.0\textwidth}
\noindent\small\textit{\bf AMS Subject Classifications:}}%
{\end{minipage}\end{center}\par}

\newenvironment{cclass}%
{\begin{center}\begin{minipage}[t]{1.0\textwidth}
\noindent\small\textit{\bf Computing Classification System codes:}}%
{\end{minipage}\end{center}\par}
\usepackage{enumerate}
\usepackage[square,sort]{natbib}

\newcommand{\MATCHtitle}%
{Adaptive Quadrilateral Mesh in Curved Domains}
\newcommand{\MATCHauthor}
{\text{Sanjay Kumar Khattri}
\\[1ex]
\text{Department of Mathematics, University of Bergen, Norway}\\
\tt{sanjay@mi.uib.no}\\ \texttt{{www.mi.uib.no/$\sim$sanjay}}}
\begin{document}
\begin{spacing}{1.0}
\maketitle
\begin{abstract}
Nonlinear elliptic system for generating adaptive quadrilateral meshes in curved domains is presented. Presented technique has been implemented in the C$^{++}$ language. The included software package can write the converged meshes in the GMV and Matlab formats. Since, grid adaptation is required for numerically capturing important characteristics of a process such as boundary layers. So, the presented technique and the software package can be a useful tool. 
\end{abstract}
\begin{class}
65M50, 65M20, 35K05
\end{class}
\begin{cclass}
G.1.0, G.4 (Mathematical Software), G.1.7 (Numerical Analysis)
\end{cclass}
\section{Introduction}
Quadrilateral grids are extensively used for numerical simulation. Accuracy of a simulation is strongly depend on the grid quality. Here, quality means orthogonality at the boundaries and quasi-orthogonality within the critical regions, smoothness, bounded aspect ratios and solution adaptive behaviour. Grid adaptation is used for increasing the efficiency of numerical schemes by focusing the computational effort where it is needed. In this article, we review the elliptic grid generation system for generating adaptive quadrilateral meshes. The presented scheme has been implemented in the C$^{++}$ language.

\begin{figure}[h!]
\begin{center}
\includegraphics[scale=0.8]{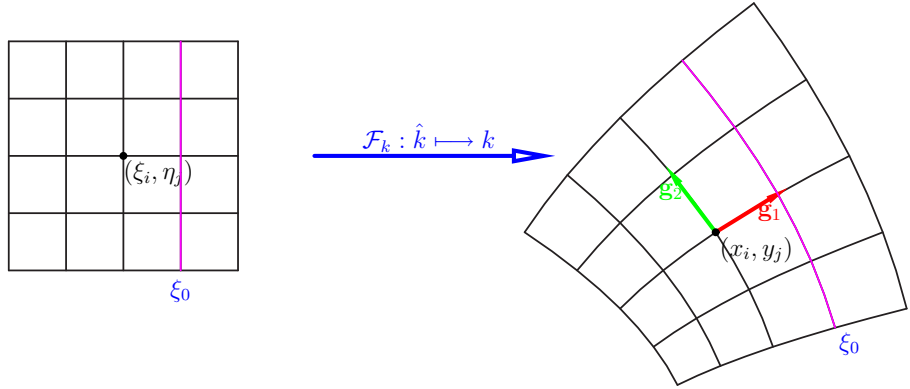}
\end{center}
\caption{Mapping $\mathcal{F}_{k}$ from a reference unit square ($\hat{k}$) on the left to a physical domain (k).} \label{map1}
\end{figure}

\begin{figure}
\begin{center}
\includegraphics[scale=0.60]{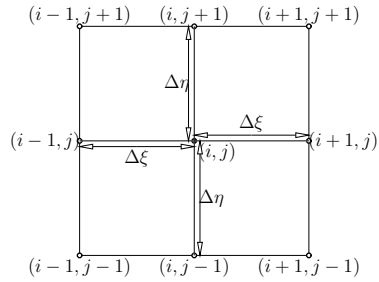}
\end{center}
\caption{Finite difference stencil in the $\xi$-$\eta$ computational space.} 
\label{fig:stencil}
\end{figure}

\begin{figure}
\begin{center}
\includegraphics[scale=1.0]{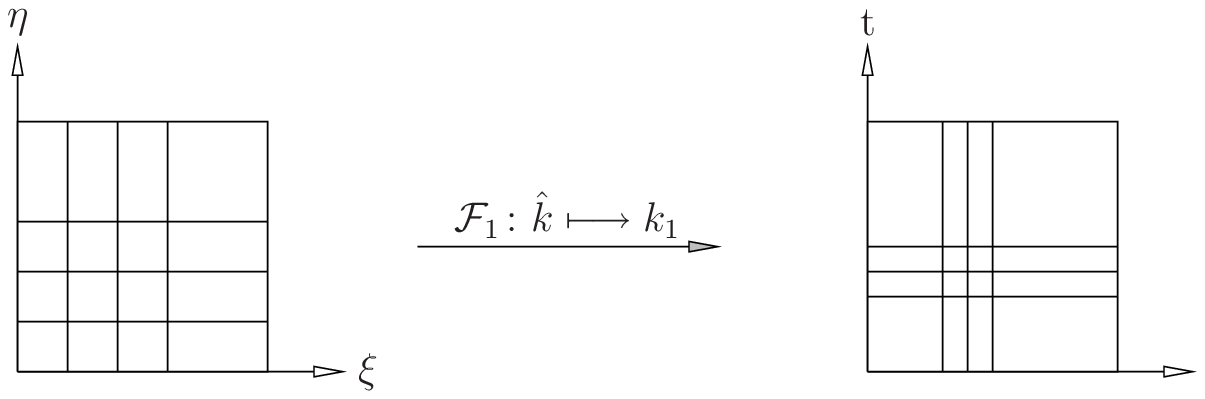}
\end{center}
\caption{Mapping $\mathcal{F}_{1}$ from a unit square ($\hat{k}$) in the reference space to a unit square in the parameter space $(k_1)$.}
\label{map2}
\end{figure}

For meshing a domain into non-simplex elements (quadrilaterals in 2D and hexahedrals in 3D), we seek a mapping from a reference square or cube to the physical domain. This mapping can be algebraic in nature such as Transfinite Interpolation or it can be expressed by a system of nonlinear partial differential equations \cite[and references therein]{knupp1,winslow} such as elliptic system. We are looking for a vector mapping, $\mathcal{F}_k(\hat{k}) = (x,y)^t$, from a unit square in the reference space ($\hat{k}=[0,1]\times[0,1]$) to a physical space ($k$); i.e. $\mathcal{F}_k\colon\hat{k}\longmapsto k$ (see Figure \ref{map1}). Mapping $\mathcal{F}_k$ gives the position of a point in the physical space corresponding to a point in the computational or reference space. Let the physical space be given by the $x$ and $y$ coordinates and the computational space be given by the $\xi$ and $\eta$ coordinates ($\xi\in[0,1]$ and $\eta\in[0,1]$). We are using the following elliptic system for defining the mapping $\mathcal{F}_k = (x,y)^t$
\begin{alignat}{2}
g_{22}\,\dfrac{\partial^2{x}}{\partial{\xi}^2}-2\,g_{12}\,\dfrac{\partial^2{x}}{\partial{\xi}{\partial\eta}}+g_{11}\,\dfrac{\partial^2{x}}{\partial\eta^2}+P\,x_{\xi}
+Q \,x_{\eta} &= 0\mathpunct{,}
\label{elliptic1}\\
g_{22}\,\dfrac{\partial^2{y}}{\partial{\xi^2}}-2\,g_{12}\,\dfrac{\partial^2{y}}{\partial{\xi}{\partial\eta}}+g_{11}\,\dfrac{\partial^2{y}}{\partial\eta^2}+P\,y_{\xi} +Q\, y_{\eta} &= 0\mathpunct{.} 
\label{elliptic2} 
\end{alignat}
Here, the terms $P$ and $Q$ are used for grid adaptation and are given as
\begin{alignat}{2}
P &= g_{22}\,P_{11}^1-2\,g_{12}\,P_{12}^1+g_{11}\,P_{22}^1, \label{P_}\\
Q &= g_{22}\,P_{11}^2-2\,g_{12}\,P_{12}^2+g_{11}\,P_{22}^2.\label{Q_}
\end{alignat}
Equations \ref{elliptic1}-\ref{elliptic2} are non-linear and are coupled through the metric coefficients $g_{ij}$ (coefficients of the metric tensor). Metric coefficients are given as 
\begin{equation}
g_{11} = x_{\xi}^2+y_{\xi}^2, \quad  g_{22} = x_{\eta}^2+y_{\eta}^2 \quad \text{and} \quad g_{12} = x_{\xi}\,x_\eta+y_{\xi}\,y_\eta. 
\end{equation}
For generating grids in the physical space, the elliptic system \ref{elliptic1}-\ref{elliptic2} is solved for the coordinates $(x,y)$ on a unit square in the computational space by the method of Finite Differences. Boundary of the physical domain is specified as the Dirichlet boundary condition on the unit square in the computational space. In the Figure \ref{map1}, $\mathbf{g}_1$ ($=\mathbf{r}_\xi$) and $\mathbf{g}_2$ ($=\mathbf{r}_\eta$) are the covariant base vectors at the point $(x_i,y_j)$. Figure \ref{fig:stencil} shows a finite difference stencil around the point ($\xi_i,\eta_j$) in the computational space. A finite difference approximation of $x_{\xi}$ and $x_{\eta}$ at the point ($i,j$) (see Figure \ref{fig:stencil}) is
$$x_{\xi}=\dfrac{\left[x(i+1,j)-x(i-1,j)\right]}{2\,\Delta{\xi}}\quad\text{and}\quad x_{\eta}=\dfrac{\left[x(i,j+1)-x(i,j-1)\right]}{2\,\Delta{\eta}}\mathpunct{.}$$
Similarly, $x_\eta$ and $y_\eta$ can be defined. 
Here, we are assuming that the grid in the computational space is uniform. However, grid in the physical space can be compressed or stretched.
Terms $P_{ij}^k$ ($i$ = 1,2 and $j$ = 1,2 and $k$ = 1,2 and $P_{12}^k = P_{21}^k$) in the equations \eqref{elliptic1}-\eqref{elliptic2} are determined through another mapping $\mathcal{F}_1$. The mapping $\mathcal{F}_1$ is shown in the Figure \ref{map2}. This mapping maps a unit square in the computational space to a unit square in the parameter space. For defining the mapping $\mathcal{F}_1\colon\hat{k}\longrightarrow{k_1}$, the boundary and internal grid points of the parameter space are mapped to the reference space. 
The Jacobian matrix $\mathbf{T}$ of the mapping $\mathcal{F}_1$ and the vectors $\mathbf{P}_{11}$, $\mathbf{P}_{12}$  and $\mathbf{P}_{22}$ are given as follows
\begin{alignat}{2}
{\mathsf{\mathbf{T}}} &=\left(\begin{array}{cc}
s_\xi & s_\eta \\
t_\xi & t_\eta 
\end{array}
\right), & \qquad \mathbf{P}_{11} &= -\mathbf{T}^{-1}\left(\begin{array}{c}
s_{\xi\xi} \\ 
t_{\xi\xi}
\end{array}
\right), \\
\mathbf{P}_{22} &= -\mathbf{T}^{-1} \left(\begin{array}{c}
s_{\eta\eta} \\ 
t_{\xi\eta}
\end{array}
\right), & \qquad \mathbf{P}_{12} &= -\mathbf{T}^{-1} \left(\begin{array}{c}
s_{\xi\eta} \\ 
t_{\xi\eta}
\end{array}
\right).
\end{alignat}
The terms $P_{ij}^1$ ($i,j = 1, 2$) are the first component of the vector $\mathbf{P}_{ij}$ and the terms $P_{ij}^2$ are the second component of the vector $\mathbf{P}_{ij}$. It should be noted that the vectors $\mathbf{P}_{11}$, $\mathbf{P}_{12}$  and $\mathbf{P}_{22}$ can be computed a priori for clustering the grid points in the physical space. A second order finite difference approximation of different operators required for computing the vectors $\mathbf{P}_{11}$, $\mathbf{P}_{22}$, $\mathbf{P}_{12}$ and the Jacobian $\mathbf{T}$ are given in the Table \ref{table1}. We are using the stencil shown in the Figure \ref{fig:stencil}.
\begin{table}
\caption{Finite difference approximation of continuous operators.}
\begin{tabular}{||c||}
\hline\\
$s_{\xi}=\dfrac{s(i+1,j)-s(i-1,j)}{2\,\Delta{\xi}}$ $,\quad$   $s_{\xi\xi}=\dfrac{s(i+1,j)-2\,s(i,j)+s(i-1,j)}
{{\Delta{\xi}}^{2}}$\\ 
\\
$t_{\xi}=\dfrac{t(i,j+1)-t(i,j-1)}{2\,\Delta{\eta}}$ $,\quad$ $t_{\eta\eta}=\dfrac{t(i,j+1)-2\,t(i,j)+t(i,j-1)}
{{\Delta{\eta}}^{2}}$\\ 
\\
$s_{\eta\eta}=\dfrac{s(i,j+1)-2\,s(i,j)+s(i,j-1)}{{\Delta{\eta}}^2}$ $,\quad$ $t_{\xi\xi}=\dfrac{t(i+1,j)-2\,t(i,j)+t(i-1,j)}{{\Delta{\xi}}^2}$\\
\\
$s_{\xi\eta}=\dfrac{s(i+1,j+1)+s(i-1,j-1)-s(i-1,j+1)-s(i+1,j-1)}{4\,\Delta{\xi}\Delta{\eta}}$\\
\\
$t_{\xi\eta}=\dfrac{t(i+1,j+1)+t(i-1,j-1)-t(i-1,j+1)-t(i+1,j-1)}{4\,\Delta{\xi}\Delta{\eta}}$\\ \\
\hline \hline
\end{tabular}
\label{table1}
\end{table}
\section{C$^{++}$ Implementation}
We have implemented the presented technique in the C$^{++}$ language for generating adaptive grids. The package can write meshes in the Matlab and GMV \cite{gmv} formats. It consists of one Domain class (see the subsections \ref{subsec_domain.h} and \ref{subsec_domain.cpp}). Domain class is used for expressing unit square in the computational space (see line no. 045 in the subsection \ref{subsec_main.cpp}), unit square in the parameter space (see line no. 022 in the subsection \ref{subsec_main.cpp}) and the physical domain (see line no. 038 in the subsection \ref{subsec_main.cpp}). The physical domain is defined in the file functions.h (see the subsection \ref{subsec_fun}). For clustering grids in the parameter space different functions are defined in the domain class (see the line numbers 026, 027, 029, 030, 032, 033 in the subsection \ref{subsec_domain.h}). 

The coupled elliptic system are linearised by the method of Finite Difference and the resulting system is solved by the SOR relaxation (see the subsection \ref{subsec_SOR}). The SOR algorithm consists of three loops, while outer loop (see line number 043 in the subsection \ref{subsec_SOR}) and two inner for loops (see line numbers 046 and 047 in the subsection \ref{subsec_SOR}). Each iteration of an inner loop provides a new mesh by the SOR relaxation. The outer loop is controlled by the maximum number of SOR iterations (see line number 028 in the subsection \ref{subsec_SOR}) and a given tolerance (see the line number 027 in the subsection \ref{subsec_SOR}). 

The overall algorithm proceeds as follows. Generate grids in the computational and parameter spaces. Compute the matrix $\mathbf{T}$ and vectors $\mathbf{P}_{ij}$ for defining the mapping $\mathcal{F}_1$ (from computational space to parameter space). An initial grid, $\mathbf{r}_{\text{old}}$, in the physical region is generated (say by Transfinite Interpolation). This information is then passed to the SOR solver (see the line number 043 in the subsection \ref{subsec_main.cpp}).
\subsection{main.cpp}
\label{subsec_main.cpp}
{\ttfamily \raggedright \tiny
001 \textsl{//+++++++++++++++++}\\
002 \#include\ <{}iostream>{}\\
003 \#include\ <{}iomanip>{}\\
004 \#include\ <{}vector>{}\\
005 \#include\ <{}iterator>{}\\
006 \#include\ <{}fstream>{}\\
007 \#include\ <{}sstream>{}\\
008 \#include\ <{}map>{}\\
009 \#include\ "{}domain.h"{}\\
010 \textsl{//\#include\ "{}write\underline\ matlab.h"{}}\\
011 \#include\ "{}matrix.h"{}\\
012 \#include\ "{}sor\underline\ solver.cpp"{}\\
013 \textsl{//+++++++++++++++++}\\
014 \textbf{int}\ main()\{\\
015 \ \ \ \ \textbf{bool}\ grid\underline\ dist\ =\ \textbf{true};\\
016 \ \ \ \ \textbf{bool}\ run\underline\ ellip\ =\ \textbf{true};\\
017 \ \ \ \ \textbf{unsigned}\ xdim,ydim\ ;\\
018 \ \ \ \ xdim\ =\ 31,ydim\ =\ 31;\\
019 \ \ \ \ \textbf{double}\ del\underline\ xi\ =\ 1.0/double(xdim-{}1.0);\\
020 \ \ \ \ \textbf{double}\ del\underline\ eta\ =\ 1.0/double(ydim-{}1.0);\\
021 \ \ \ \ \textsl{//Parameter\ Space\ ref(xdim,ydim);}\\
022 \ \ \ \ Domain\ parm(xdim,ydim);\\
023 \ \ \ \ \textsl{//Meshing\ the\ Parameter\ }\\
024 \ \ \ \ parm.Grid\underline\ Gen();\\
025 \ \ \ \ \textsl{//Clustering\ the\ Mesh}\\
026 \ \ \ \ \textsl{//Example\ 1}\\
027 \ \ \ \ \textsl{//parm.Cluster\underline\ X\underline\ Near(0.5);}\\
028 \ \ \ \ \textsl{//parm.Cluster\underline\ Y\underline\ Near(0.5);}\\
029 \ \ \ \ \textsl{//Example\ 2}\\
030 \ \ \ \ \textsl{//parm.Cluster\underline\ Two\underline\ Lines\underline\ X(0.25,0.750);}\\
031 \ \ \ \ \textsl{//parm.Cluster\underline\ Two\underline\ Lines\underline\ Y(0.25,0.750);}\\
032 \ \ \ \ \textsl{//Example\ 3}\\
033 \ \ \ \ parm.Bound\underline\ Clust\underline\ X(0.5);\\
034 \ \ \ \ parm.Bound\underline\ Clust\underline\ Y(0.5);\\
035 \ \ \ \ \ \ \ \ \\
036 \ \ \ \ parm.Fill\underline\ del\underline\ xi\underline\ eta(del\underline\ xi,del\underline\ eta);\\
037 \ \ \ \ \ \ \ \ \\
038 \ \ \ \ Domain\ physical(xdim,ydim);\\
039 \ \ \ \ physical.Read\underline\ Bd();\\
040 \ \ \ \ physical.Fill\underline\ del\underline\ xi\underline\ eta(del\underline\ xi,del\underline\ eta);\\
041 \ \ \ \ \textbf{unsigned}\ max\underline\ iter\ =\ 100;\\
042 \ \ \ \ \textbf{double}\ w\ =\ 1.90;\\
043 \ \ \ \ SORSOLVER(physical,\ parm,\ xdim,\ ydim);\\
044 \ \ \ \ \textsl{//Reference\ or\ computational\ space}\\
045 \ \ \ \ Domain\ ref(xdim,ydim);\\
046 \ \ \ \ \\
047 \ \ \ \ \textsl{//Writing\ the\ mesh\ in\ the\ physical\ space\ (GMV)}\\
048 \ \ \ \ std::ofstream\ outPhy("{}gmv\underline\ Physical.dat"{},std::ios::out);\\
049 \ \ \ \ \textbf{if}(!(outPhy))\ std::cerr\ <{}<{}\ "{}ERROR\ :\ UNABLE\ TO\ OPEN\ $\backslash$''outPhy$\backslash$''$\backslash$n"{};\\
050 \ \ \ \ physical.GMV\underline\ Writer(outPhy);\\
051 \ \ \ \ \textbf{if}(outPhy.is\underline\ open())\ \ outPhy.close();\\
052 \ \ \ \ \\
053 \ \ \ \ \textsl{//writing\ mesh\ in\ the\ parameter\ space\ (GMV)}\\
054 \ \ \ \ std::ofstream\ outParm("{}gmv\underline\ Para.dat"{},std::ios::out);\\
055 \ \ \ \ \textbf{if}(!(outParm))\ std::cerr\ <{}<{}\ "{}ERROR\ :\ UNABLE\ TO\ OPEN\ $\backslash$''outParm$\backslash$''$\backslash$n"{};\\
056 \ \ \ \ parm.GMV\underline\ Writer(outParm);\\
057 \ \ \ \ \textbf{if}(outParm.is\underline\ open())\ \ outParm.close();\\
058 \ \ \ \ parm.Matlab\underline\ Writer();\\
059 \ \ \ \ \\
060 \ \ \ \ \textbf{return}\ EXIT\underline\ SUCCESS;\\
061 \}\\
062 \ \\
063  }
\normalfont\normalsize

\subsection{domain.h}
\label{subsec_domain.h}
{\ttfamily \raggedright \tiny
001 \#ifndef PARAMETER\underline\ SPACE\\
002 \#define PARAMETER\underline\ SPACE\\
003 \textsl{//+++++++++++++++++}\\
004 \#include<{}iostream>{}\\
005 \#include<{}iomanip>{}\\
006 \#include<{}vector>{}\\
007 \#include<{}iterator>{}\\
008 \#include<{}fstream>{}\\
009 \#include<{}sstream>{}\\
010 \#include<{}map>{}\\
011 \textsl{//++++++++++++++++}\\
012 \#include\ "{}matrix.h"{}\\
013 \textsl{//+++++++++++++++}\\
014 \textbf{class}\ Domain\{\\
015 \textbf{public}:\\
016 \ \ \ \ Domain();\\
017 \ \ \ \ Domain(\textbf{unsigned}\ \textbf{int}\ xdim1,\ \textbf{unsigned}\ ydim1);\\
018 \ \ \ \ Domain(\textbf{const}\ Domain\ \&\ org);\\
019 \ \ \ \ \textbf{unsigned}\ \textbf{int}\ XDIM()\ \textbf{const}\ ;\\
020 \ \ \ \ \textbf{unsigned}\ \textbf{int}\ YDIM()\ \textbf{const}\ ;\\
021 \ \ \ \ \textbf{void}\ Grid\underline\ Gen();\\
022 \ \ \ \ std::vector<{}\textbf{double}>{}\ XCOORDS();\\
023 \ \ \ \ std::vector<{}\textbf{double}>{}\ YCOORDS();\\
024 \ \ \ \ \textbf{double}\ Eriksson\underline\ 1(\textbf{double}\ eta);\\
025 \ \ \ \ \\
026 \ \ \ \ \textbf{void}\ Cluster\underline\ X\underline\ Near(\textbf{double}\ eta0);\\
027 \ \ \ \ \textbf{void}\ Cluster\underline\ Y\underline\ Near(\textbf{double}\ eta0);\\
028 \ \ \ \ \\
029 \ \ \ \ \textbf{void}\ Cluster\underline\ Two\underline\ Lines\underline\ X(\textbf{double}\ eta1,\textbf{double}\ eta2);\\
030 \ \ \ \ \textbf{void}\ Cluster\underline\ Two\underline\ Lines\underline\ Y(\textbf{double}\ eta1,\textbf{double}\ eta2);\\
031 \ \ \ \ \\
032 \ \ \ \ \textbf{void}\ Bound\underline\ Clust\underline\ X(\textbf{double}\ eta1);\\
033 \ \ \ \ \textbf{void}\ Bound\underline\ Clust\underline\ Y(\textbf{double}\ eta1);\\
034 \ \ \ \ \\
035 \ \ \ \ \textbf{double}\&\ XCOORD(\textbf{unsigned}\ \textbf{int}\ i\ ,\ \textbf{unsigned}\ j);\\
036 \ \ \ \ \textbf{double}\&\ YCOORD(\textbf{unsigned}\ \textbf{int}\ i\ ,\ \textbf{unsigned}\ j);\\
037 \ \ \ \ \textbf{void}\ Read\underline\ Bd();\\
038 \ \ \ \ \textbf{void}\ Matlab\underline\ Writer();\\
039 \ \ \ \ \textbf{void}\ GMV\underline\ Writer(std::ofstream\ \&\ outFile);\\
040 \ \ \ \ \textbf{void}\ Fill\underline\ del\underline\ xi\underline\ eta(\textbf{double}\ xi,\textbf{double}\ eta);\\
041 \ \ \ \ \textsl{//void\ Call\underline\ Grid\underline\ Adapter();}\\
042 \ \ \ \ \\
043 \ \ \ \ \textsl{//+++++++++++++}\\
044 \ \ \ \ std::vector<{}\textbf{double}>{}\ P11(\textbf{unsigned}\ \textbf{int}\ i\ ,\ \textbf{unsigned}\ \textbf{int}\ j);\\
045 \ \ \ \ std::vector<{}\textbf{double}>{}\ P22(\textbf{unsigned}\ \textbf{int}\ i\ ,\ \textbf{unsigned}\ \textbf{int}\ j);\\
046 \ \ \ \ std::vector<{}\textbf{double}>{}\ P12(\textbf{unsigned}\ \textbf{int}\ i\ ,\ \textbf{unsigned}\ \textbf{int}\ j);\\
047 \ \ \ \ \textsl{//++++++++++}\\
048 \ \ \ \ Matrix\ MeshX();\\
049 \ \ \ \ Matrix\ MeshY();\\
050 \ \ \ \ \textsl{//+++++++++++++}\\
051 \ \ \ \ \textbf{double}\ G22(\textbf{unsigned}\ \textbf{int}\ i\ ,\ \textbf{unsigned}\ \textbf{int}\ j);\\
052 \ \ \ \ \textbf{double}\ G11(\textbf{unsigned}\ \textbf{int}\ i\ ,\ \textbf{unsigned}\ \textbf{int}\ j);\\
053 \ \ \ \ \textbf{double}\ X\underline\ xi(\textbf{unsigned}\ \textbf{int}\ i\ ,\ \textbf{unsigned}\ \textbf{int}\ j);\\
054 \ \ \ \ \textbf{double}\ Y\underline\ xi(\textbf{unsigned}\ \textbf{int}\ i\ ,\ \textbf{unsigned}\ \textbf{int}\ j);\\
055 \ \ \ \ \textbf{double}\ X\underline\ eta(\textbf{unsigned}\ \textbf{int}\ i\ ,\ \textbf{unsigned}\ \textbf{int}\ j);\\
056 \ \ \ \ \textbf{double}\ Y\underline\ eta(\textbf{unsigned}\ \textbf{int}\ i\ ,\ \textbf{unsigned}\ \textbf{int}\ j);\\
057 \ \ \ \ \textbf{double}\ X\underline\ xieta(\textbf{unsigned}\ \textbf{int}\ i\ ,\ \textbf{unsigned}\ \textbf{int}\ j);\\
058 \ \ \ \ \textbf{double}\ Y\underline\ xieta(\textbf{unsigned}\ \textbf{int}\ i\ ,\ \textbf{unsigned}\ \textbf{int}\ j);\\
059 \ \ \ \ \\
060 \textbf{private}:\\
061 \ \ \ \ \textbf{double}\ del\underline\ eta,del\underline\ xi;\\
062 \ \ \ \ \textbf{unsigned}\ xdim,ydim;\\
063 \ \ \ \ Matrix\ x,y;\\
064 \ \ \ \ std::vector<{}\textbf{double}>{}\ xcoords,ycoords;\\
065 \};\\
066 \#endif\\
067 \ \\
068  }
\normalfont\normalsize

\subsection{domain.cpp}
\label{subsec_domain.cpp}
{\ttfamily \raggedright \tiny
001 \#include\ "{}domain.h"{}\\
002 \ \\
003 \#ifndef \underline\ FUNCTIONS\underline\ \\
004 \#include\ "{}functions.h"{}\ \\
005 \#endif\\
006 \ \\
007 \#include\ <{}cassert>{}\\
008 \ \\
009 Domain::Domain()\{\\
010 \ \ \ \ xdim\ =\ 0\ ;\ ydim\ =\ 0;\\
011 \}\\
012 Domain::Domain(\textbf{unsigned}\ \textbf{int}\ xdim1,\textbf{unsigned}\ \textbf{int}\ ydim1)\{\\
013 \ \ \ \ xdim\ =\ xdim1\ ;\ ydim\ =\ ydim1;\\
014 \}\\
015 Domain::Domain(\textbf{const}\ Domain\ \&\ org)\{\\
016 \ \ \ \ xdim\ =\ org.XDIM();\\
017 \ \ \ \ ydim\ =\ org.YDIM();\\
018 \ \ \ \ Grid\underline\ Gen();\\
019 \}\\
020 \textbf{unsigned}\ \textbf{int}\ Domain::XDIM()\ \textbf{const}\{\\
021 \ \ \ \ \textbf{return}\ xdim;\\
022 \}\\
023 \textbf{unsigned}\ \textbf{int}\ Domain::YDIM()\ \textbf{const}\{\\
024 \ \ \ \ \textbf{return}\ ydim;\\
025 \}\\
026 \textbf{void}\ Domain::Grid\underline\ Gen()\{\\
027 \ \ \ \ Matrix\ xt(xdim,ydim),yt(xdim,ydim);\\
028 \ \ \ \ assert(0\ !=\ xdim\ \&\&\ 0\ !=\ ydim);\\
029 \ \ \ \ \textbf{for}(\textbf{unsigned}\ \textbf{int}\ j\ =\ 0\ ;\ j\ <{}\ ydim\ ;\ ++j)\{\\
030 \ \ \textbf{for}(\textbf{unsigned}\ \textbf{int}\ i\ =\ 0\ ;\ i\ <{}\ xdim\ ;\ ++i)\{\\
031 \ \ \ \ \ \ \textbf{double}\ t\underline\ x\ =\ \textbf{double}(i)/double(xdim-{}1.0);\\
032 \ \ \ \ \ \ \textbf{double}\ t\underline\ y\ \ =\ \textbf{double}(j)/double(ydim-{}1.0);\\
033 \ \ \ \ \ \ xt(i,j)\ =\ t\underline\ x\ ;\ yt(i,j)\ =\ t\underline\ y;\\
034 \ \ \}\\
035 \ \ \ \ \}\\
036 \ \ \ \ x\ =\ xt\ ;\ y\ =\ yt;\\
037 \}\\
038 std::vector<{}\textbf{double}>{}\ Domain::XCOORDS()\{\\
039 \ \ \ \ xcoords.resize(xdim$\ast$ydim);\ ycoords.resize(xdim$\ast$ydim);\\
040 \ \ \ \ \textbf{for}(\textbf{int}\ j\ =\ 0\ ;\ j\ <{}\ ydim\ ;\ ++j)\{\\
041 \ \ \textbf{for}(\textbf{int}\ i\ =\ 0\ ;\ i\ <{}\ xdim\ ;\ ++i)\{\\
042 \ \ \ \ \ \ \textbf{int}\ no\ =\ i+j$\ast$xdim;\\
043 \ \ \ \ \ \ xcoords[no]\ =\ x(i,j);\\
044 \ \ \ \ \ \ ycoords[no]\ =\ y(i,j);\\
045 \ \ \}\\
046 \ \ \ \ \}\\
047 \ \ \ \ \textbf{return}\ xcoords;\\
048 \}\\
049 std::vector<{}\textbf{double}>{}\ Domain::YCOORDS()\{\\
050 \ \ \ \ xcoords.resize(xdim$\ast$ydim);\ ycoords.resize(xdim$\ast$ydim);\\
051 \ \ \ \ \textbf{for}(\textbf{int}\ j\ =\ 0\ ;\ j\ <{}\ ydim\ ;\ ++j)\{\\
052 \ \ \textbf{for}(\textbf{int}\ i\ =\ 0\ ;\ i\ <{}\ xdim\ ;\ ++i)\{\\
053 \ \ \ \ \ \ \textbf{int}\ no\ =\ i+j$\ast$xdim;\\
054 \ \ \ \ \ \ xcoords[no]\ =\ x(i,j);\\
055 \ \ \ \ \ \ ycoords[no]\ =\ y(i,j);\\
056 \ \ \}\\
057 \ \ \ \ \}\\
058 \ \ \ \ \textbf{return}\ ycoords;\\
059 \}\\
060 \ \\
061 \textbf{void}\ Domain::Bound\underline\ Clust\underline\ X(\textbf{double}\ eta1)\{\\
062 \ \\
063 \ \ \ \ \textbf{double}\ alpha\ =\ 4.0;\\
064 \ \ \ \ \textbf{double}\ h\ =\ 1.0;\\
065 \ \ \ \ \textbf{double}\ h2\ =\ 1.0;\\
066 \ \ \ \ \textbf{double}\ h1\ =\ 0.0;\\
067 \ \ \\
068 \ \ \ \ \textbf{for}(\textbf{int}\ j\ \ =\ 0\ ;\ j\ <{}\ ydim\ ;\ ++j)\{\\
069 \ \ \textbf{for}(\textbf{int}\ i\ =\ 0\ ;\ i\ <{}\ xdim\ ;\ ++i)\{\\
070 \ \ \ \ \ \ \\
071 \ \ \ \ \ \ \textbf{if}(x(i,j)\ <{}=\ eta1\ \&\&\ 0\ <{}=\ x(i,j))\{\\
072 \ \ \ \ \textbf{double}\ eta\ =\ x(i,j);\\
073 \ \ \ \ x(i,j)\ =(h2-{}h1)$\ast$\ eta1$\ast$(std::exp(alpha$\ast$eta/eta1)-{}1.0)/(std::exp(alpha)-{}1.0)+h1;\ \ \\
074 \ \ \ \ \ \ \}\\
075 \ \ \ \ \ \ \\
076 \ \ \ \ \ \ \textbf{if}(x(i,j)\ >{}=\ eta1\ \&\&\ x(i,j)\ <{}=\ 1.0)\{\\
077 \ \ \ \ \textbf{double}\ eta\ =\ x(i,j);\\
078 \ \ \ \ x(i,j)\ =\ (h2-{}h1)$\ast$(1.0-{}(1.0-{}eta1)$\ast$(((std::exp(alpha$\ast$(1.0-{}eta)/(1.0-{}eta1)))-{}1.0)/(std::exp(alpha)-{}1.0)));\\
079 \ \ \ \ \ \ \}\\
080 \ \ \ \ \ \ \\
081 \ \ \}\\
082 \ \ \ \ \}\\
083 \ \ \ \ \\
084 \}\\
085 \ \\
086 \textbf{void}\ Domain::Bound\underline\ Clust\underline\ Y(\textbf{double}\ eta1)\{\\
087 \ \\
088 \ \ \ \ \textbf{double}\ alpha\ =\ 4.0;\\
089 \ \ \ \ \textbf{double}\ h\ =\ 1.0;\\
090 \ \ \ \ \textbf{double}\ h2\ =\ 1.0;\\
091 \ \ \ \ \textbf{double}\ h1\ =\ 0.0;\\
092 \ \ \\
093 \ \ \ \ \textbf{for}(\textbf{int}\ j\ \ =\ 0\ ;\ j\ <{}\ ydim\ ;\ ++j)\{\\
094 \ \ \textbf{for}(\textbf{int}\ i\ =\ 0\ ;\ i\ <{}\ xdim\ ;\ ++i)\{\\
095 \ \ \ \ \ \ \\
096 \ \ \ \ \ \ \textbf{if}(y(i,j)\ <{}=\ eta1\ \&\&\ 0\ <{}=\ y(i,j))\{\\
097 \ \ \ \ \textbf{double}\ eta\ =\ y(i,j);\\
098 \ \ \ \ y(i,j)\ =(h2-{}h1)$\ast$\ eta1$\ast$(std::exp(alpha$\ast$eta/eta1)-{}1.0)/(std::exp(alpha)-{}1.0)+h1;\ \ \\
099 \ \ \ \ \ \ \}\\
100 \ \ \ \ \ \ \\
101 \ \ \ \ \ \ \textbf{if}(y(i,j)\ >{}=\ eta1\ \&\&\ y(i,j)\ <{}=\ 1.0)\{\\
102 \ \ \ \ \textbf{double}\ eta\ =\ y(i,j);\\
103 \ \ \ \ y(i,j)\ =\ (h2-{}h1)$\ast$(1.0-{}(1.0-{}eta1)$\ast$(((std::exp(alpha$\ast$(1.0-{}eta)/(1.0-{}eta1)))-{}1.0)/(std::exp(alpha)-{}1.0)));\\
104 \ \ \ \ \ \ \}\\
105 \ \ \ \ \ \ \\
106 \ \ \}\\
107 \ \ \ \ \}\\
108 \ \ \ \ \\
109 \}\\
110 \ \\
111 \textsl{//===========}\\
112 \textbf{double}\ Domain::Eriksson\underline\ 1(\textbf{double}\ eta)\{\\
113 \ \ \ \ \textbf{double}\ h\ =\ 1.0;\\
114 \ \ \ \ \textbf{double}\ alpha\ =\ 3.0;\\
115 \ \\
116 \ \ \ \ \textbf{return}\ h$\ast$((std::exp(alpha$\ast$eta)-{}1.0)/(std::exp(alpha)-{}1.0));\\
117 \ \ \ \ \\
118 \ \ \ \ \\
119 \}\\
120 \ \\
121 \textbf{void}\ Domain::Cluster\underline\ Two\underline\ Lines\underline\ X(\textbf{double}\ eta1,\textbf{double}\ eta2)\{\\
122 \ \\
123 \ \ \ \ \textbf{double}\ alpha\ =\ 5.0;\\
124 \ \ \ \ \textbf{double}\ h\ =\ 1.0;\\
125 \ \ \ \ \textbf{double}\ eta0\ =\ (eta1+eta2)$\ast$0.5;\\
126 \ \ \\
127 \ \ \ \ \textbf{for}(\textbf{int}\ j\ \ =\ 0\ ;\ j\ <{}\ ydim\ ;\ ++j)\{\\
128 \ \ \textbf{for}(\textbf{int}\ i\ =\ 0\ ;\ i\ <{}\ xdim\ ;\ ++i)\{\\
129 \ \ \ \ \ \ \\
130 \ \ \ \ \ \ \textbf{if}(x(i,j)\ <{}=\ eta1\ \&\&\ 0\ <{}=\ x(i,j))\{\\
131 \ \ \ \ \textbf{double}\ eta\ =\ x(i,j);\\
132 \ \ \ \ x(i,j)\ =\ eta1$\ast$(h-{}Eriksson\underline\ 1(1-{}eta/eta1));\ \ \\
133 \ \ \ \ \ \ \}\\
134 \ \ \ \ \ \ \\
135 \ \ \ \ \ \ \textbf{if}(x(i,j)\ >{}=\ eta1\ \&\&\ x(i,j)\ <{}=\ eta0\ )\{\\
136 \ \ \ \ \textbf{double}\ eta\ =\ x(i,j);\\
137 \ \ \ \ x(i,j)\ =\ h$\ast$eta1+\ (eta0-{}eta1)$\ast$Eriksson\underline\ 1((eta-{}eta1)/(eta0-{}eta1));\\
138 \ \ \ \ \\
139 \ \ \ \ \ \ \}\\
140 \ \ \ \ \ \ \\
141 \ \ \ \ \ \ \textbf{if}(x(i,j)\ >{}=\ eta0\ \&\&\ x(i,j)\ <{}=\ eta2)\{\\
142 \ \ \ \ \textbf{double}\ eta\ =\ x(i,j);\\
143 \ \ \ \ x(i,j)\ =\ h$\ast$eta0\ +\ (eta2-{}eta0)$\ast$(h-{}Eriksson\underline\ 1((eta2-{}eta)/(eta2-{}eta0)));\\
144 \ \ \ \ \ \ \}\\
145 \ \ \ \ \ \ \\
146 \ \ \ \ \ \ \textbf{if}(x(i,j)\ >{}=\ eta2\ \&\&\ x(i,j)\ <{}=\ 1.0)\{\\
147 \ \ \ \ \textbf{double}\ eta\ =\ x(i,j);\\
148 \ \ \ \ x(i,j)\ =\ h$\ast$eta2\ +\ (1.0-{}eta2)$\ast$Eriksson\underline\ 1((eta-{}eta2)/(1.0-{}eta2));\\
149 \ \ \ \ \ \ \}\\
150 \ \ \ \ \ \ \\
151 \ \ \}\\
152 \ \ \ \ \}\\
153 \ \ \ \ \\
154 \}\\
155 \ \\
156 \textbf{void}\ Domain::Cluster\underline\ Two\underline\ Lines\underline\ Y(\textbf{double}\ eta1,\textbf{double}\ eta2)\{\\
157 \ \\
158 \ \ \ \ \textbf{double}\ alpha\ =\ 5.0;\\
159 \ \ \ \ \textbf{double}\ h\ =\ 1.0;\\
160 \ \ \ \ \textbf{double}\ eta0\ =\ (eta1+eta2)$\ast$0.5;\\
161 \ \ \\
162 \ \ \ \ \textbf{for}(\textbf{int}\ j\ \ =\ 0\ ;\ j\ <{}\ ydim\ ;\ ++j)\{\\
163 \ \ \textbf{for}(\textbf{int}\ i\ =\ 0\ ;\ i\ <{}\ xdim\ ;\ ++i)\{\\
164 \ \ \ \ \ \ \\
165 \ \ \ \ \ \ \textbf{if}(y(i,j)\ <{}=\ eta1\ \&\&\ 0\ <{}=\ y(i,j))\{\\
166 \ \ \ \ \textbf{double}\ eta\ =\ y(i,j);\\
167 \ \ \ \ y(i,j)\ =\ eta1$\ast$(h-{}Eriksson\underline\ 1(1-{}eta/eta1));\ \ \\
168 \ \ \ \ \ \ \}\\
169 \ \ \ \ \ \ \\
170 \ \ \ \ \ \ \textbf{if}(y(i,j)\ >{}=\ eta1\ \&\&\ y(i,j)\ <{}=\ eta0\ )\{\\
171 \ \ \ \ \textbf{double}\ eta\ =\ y(i,j);\\
172 \ \ \ \ y(i,j)\ =\ h$\ast$eta1+\ (eta0-{}eta1)$\ast$Eriksson\underline\ 1((eta-{}eta1)/(eta0-{}eta1));\\
173 \ \ \ \ \\
174 \ \ \ \ \ \ \}\\
175 \ \ \ \ \ \ \\
176 \ \ \ \ \ \ \textbf{if}(y(i,j)\ >{}=\ eta0\ \&\&\ y(i,j)\ <{}=\ eta2)\{\\
177 \ \ \ \ \textbf{double}\ eta\ =\ y(i,j);\\
178 \ \ \ \ y(i,j)\ =\ h$\ast$eta0\ +\ (eta2-{}eta0)$\ast$(h-{}Eriksson\underline\ 1((eta2-{}eta)/(eta2-{}eta0)));\\
179 \ \ \ \ \ \ \}\\
180 \ \ \ \ \ \ \\
181 \ \ \ \ \ \ \textbf{if}(y(i,j)\ >{}=\ eta2\ \&\&\ y(i,j)\ <{}=\ 1.0)\{\\
182 \ \ \ \ \textbf{double}\ eta\ =\ y(i,j);\\
183 \ \ \ \ y(i,j)\ =\ h$\ast$eta2\ +\ (1.0-{}eta2)$\ast$Eriksson\underline\ 1((eta-{}eta2)/(1.0-{}eta2));\\
184 \ \ \ \ \ \ \}\\
185 \ \ \ \ \ \ \\
186 \ \ \}\\
187 \ \ \ \ \}\\
188 \ \ \ \ \\
189 \}\\
190 \textsl{//=============}\\
191 \textbf{void}\ Domain::Cluster\underline\ X\underline\ Near(\textbf{double}\ eta0)\{\\
192 \ \ \ \ \textbf{double}\ alpha\ =\ 3.0;\\
193 \ \ \ \ \textbf{for}(\textbf{int}\ j\ \ =\ 0\ ;\ j\ <{}\ ydim\ ;\ ++j)\{\\
194 \ \ \textbf{for}(\textbf{int}\ i\ =\ 0\ ;\ i\ <{}\ xdim\ ;\ ++i)\{\\
195 \ \ \ \ \ \ \textbf{if}(x(i,j)\ <{}\ eta0)\{\\
196 \ \ \ \ \textbf{double}\ eta\ =\ x(i,j);\\
197 \ \ \ \ x(i,j)\ =\ (\textbf{double})\ eta0\ $\ast$\ (exp(alpha)\ -{}\ exp(alpha\ $\ast$\ (\textbf{double})\ (1\ -{}\ eta\ /\ eta0)))\ /\ (exp(alpha)\ -{}\ 0.1e1);\\
198 \ \ \ \ \ \ \}\\
199 \ \ \ \ \ \ \textbf{if}(x(i,j)\ >{}\ eta0)\{\\
200 \ \ \ \ \textbf{double}\ eta\ =\ x(i,j);\\
201 \ \ \ \ x(i,j)\ =\ (\textbf{double})\ eta0\ +\ (\textbf{double})\ (1\ -{}\ eta0)\ $\ast$\ (exp((\textbf{double})\ (alpha\ $\ast$\ (eta\ -{}\ eta0)\ /\ (1\ -{}\ eta0)))\ -{}\ 0.1e1)\\
202 \ \ \ \ \ \ \ \ /\ (exp((\textbf{double})\ alpha)\ -{}\ 0.1e1);\\
203 \ \ \ \ \ \ \}\\
204 \ \ \}\\
205 \ \ \ \ \}\\
206 \ \ \ \ \\
207 \}\\
208 \ \\
209 \textbf{void}\ Domain::Cluster\underline\ Y\underline\ Near(\textbf{double}\ eta0)\{\\
210 \ \ \ \ \textbf{double}\ alpha\ =\ 3.0;\\
211 \ \ \ \ \textbf{for}(\textbf{int}\ j\ \ =\ 0\ ;\ j\ <{}\ ydim\ ;\ ++j)\{\\
212 \ \ \textbf{for}(\textbf{int}\ i\ =\ 0\ ;\ i\ <{}\ xdim\ ;\ ++i)\{\\
213 \ \ \ \ \ \ \textbf{if}(y(i,j)\ <{}\ eta0)\{\\
214 \ \ \ \ \textbf{double}\ eta\ =\ y(i,j);\\
215 \ \ \ \ y(i,j)\ =\ (\textbf{double})\ eta0\ $\ast$\ (exp(alpha)\ -{}\ exp(alpha\ $\ast$\ (\textbf{double})\ (1\ -{}\ eta\ /\ eta0)))\ /\ (exp(alpha)\ -{}\ 0.1e1);\\
216 \ \ \ \ \ \ \}\\
217 \ \ \ \ \ \ \textbf{if}(y(i,j)\ >{}\ eta0)\{\\
218 \ \ \ \ \textbf{double}\ eta\ =\ y(i,j);\\
219 \ \ \ \ y(i,j)\ =\ (\textbf{double})\ eta0\ +\ (\textbf{double})\ (1\ -{}\ eta0)\ $\ast$\ (exp((\textbf{double})\ (alpha\ $\ast$\ (eta\ -{}\ eta0)\ /\ (1\ -{}\ eta0)))\ -{}\ 0.1e1)\\
220 \ \ \ \ \ \ \ \ /\ (exp((\textbf{double})\ alpha)\ -{}\ 0.1e1);\\
221 \ \ \ \ \ \ \}\\
222 \ \ \}\\
223 \ \ \ \ \}\\
224 \}\\
225 \textbf{double}\&\ Domain::XCOORD(\textbf{unsigned}\ \textbf{int}\ i\ ,\ \textbf{unsigned}\ \textbf{int}\ j)\ \{\\
226 \ \ \ \ \textbf{if}(i\ >{}=\ xdim\ ||\ j\ >{}=\ ydim\ ||\ i\ <{}\ 0\ ||\ j\ <{}\ 0)\{\\
227 \ \ std::cerr\ <{}<{}\ "{}In\ XCOORDS(..,..)\ dim\ mismatch$\backslash$n"{};\\
228 \ \ \ \ \}\\
229 \ \ \ \ \textbf{return}\ x(i,j);\\
230 \}\\
231 \textbf{double}\&\ Domain::YCOORD(\textbf{unsigned}\ \textbf{int}\ i\ ,\ \textbf{unsigned}\ \textbf{int}\ j)\{\\
232 \ \ \ \ \textbf{if}(i\ >{}=\ xdim\ ||\ j\ >{}=\ ydim\ ||\ i\ <{}\ 0\ ||\ j\ <{}\ 0)\{\\
233 \ \ std::cerr\ <{}<{}\ "{}In\ YCOORDS(..,..)\ dim\ mismatch$\backslash$n"{};\\
234 \ \ \ \ \}\\
235 \ \ \ \ \textbf{return}\ y(i,j);\\
236 \}\\
237 \textbf{void}\ Domain::Read\underline\ Bd()\{\\
238 \textsl{//Read\ the\ boundary\ of\ the\ physical\ domain}\\
239 \ \ \ \ Matrix\ \ xt(xdim,ydim),yt(xdim,ydim);\\
240 \ \ \ \ \textbf{for}(\textbf{int}\ j\ =\ 0\ ;\ j\ <{}\ ydim\ ;\ ++j)\{\\
241 \ \ \textbf{for}(\textbf{int}\ i\ =0\ ;\ i\ <{}\ xdim\ ;\ ++i)\{\\
242 \ \ \ \ \ \ \textbf{if}(0\ ==\ i\ ||\ xdim-{}1\ ==\ i\ ||\ 0\ ==\ j\ ||\ ydim-{}1\ ==\ j)\{\\
243 \ \ \ \ \textbf{double}\ zeta1\ =\ \textbf{double}(i)/double(xdim-{}1.0);\\
244 \ \ \ \ \textbf{double}\ eta1\ \ =\ \textbf{double}(j)/double(ydim-{}1.0);\\
245 \ \ \ \ xt(i,j)\ =\ zeta1;\ yt(i,j)\ =\ eta1;\\
246 \ \ \ \ \\
247 \ \ \ \ xt(i,j)\ =\ XYcircle(zeta1,eta1,1);\\
248 \ \ \ \ yt(i,j)\ =\ XYcircle(zeta1,eta1,2);\\
249 \ \ \ \ \\
250 \ \ \ \ \ \ \}\\
251 \ \ \}\\
252 \ \ \ \ \}\\
253 \textsl{//Create\ grid\ by\ the\ TFI\ \ \ \ }\\
254 \ \ \ \ \textbf{for}(\textbf{int}\ j\ =\ 1\ ;\ j\ <{}\ ydim-{}1\ ;\ ++j)\{\\
255 \ \ \textbf{for}(\textbf{int}\ i\ =\ 1\ ;\ i\ <{}\ xdim-{}1\ ;\ ++i)\{\\
256 \ \ \ \ \ \ \textbf{double}\ zeta1\ =\ \textbf{double}(i)/double(xdim-{}1.0);\\
257 \ \ \ \ \ \ \textbf{double}\ eta1\ \ =\ \textbf{double}(j)/double(ydim-{}1.0);\\
258 \ \ \ \ \ \ xt(i,j)\ =\ (1.0-{}zeta1)$\ast$xt(0,j)\ \ +zeta1$\ast$xt(xdim-{}1,j)\ +\ (1.0-{}eta1)$\ast$xt(i,0)+eta1$\ast$xt(i,ydim-{}1)-{}\\
259 \ \ \ \ ((1.0-{}zeta1)$\ast$(1.0-{}eta1)$\ast$xt(0,0)\ \ \ \ \ \ \ \ \ \ \ \ \ \ +\ (zeta1)$\ast$(1.0-{}eta1)$\ast$xt(xdim-{}1,0)\\
260 \ \ \ \ \ +\ (zeta1)$\ast$(eta1)$\ast$xt(xdim-{}1,ydim-{}1)\ \ \ \ \ \ \ \ \ \ +\ (1.0-{}zeta1)$\ast$eta1$\ast$xt(0,ydim-{}1));\\
261 \ \ \ \ \ \ yt(i,j)\ =\ (1-{}zeta1)$\ast$yt(0,j)+zeta1$\ast$yt(xdim-{}1,j)\ +\ (1.0-{}eta1)$\ast$yt(i,0)+eta1$\ast$yt(i,ydim-{}1)-{}\\
262 \ \ \ \ ((1-{}zeta1)$\ast$(1-{}eta1)$\ast$yt(0,0)\ \ \ \ \ \ \ \ \ \ \ \ \ \ \ +\ (zeta1)$\ast$(1-{}eta1)$\ast$yt(xdim-{}1,0)\\
263 \ \ \ \ \ +\ (zeta1)$\ast$(eta1)$\ast$yt(xdim-{}1,ydim-{}1)\ \ \ \ \ \ \ +\ (1-{}zeta1)$\ast$eta1$\ast$yt(0,ydim-{}1));\\
264 \ \ \}\\
265 \ \ \ \ \}\\
266 \ \ \ \ x\ =\ xt\ ;\ y\ =\ yt\ ;\ \\
267 \}\\
268 \textsl{//=============}\\
269 \textbf{void}\ Domain::Matlab\underline\ Writer()\{\\
270 \ \ \ \ \\
271 \ \ \ \ std::vector<{}\textbf{double}>{}\ x1\ =\ XCOORDS();\\
272 \ \ \ \ std::vector<{}\textbf{double}>{}\ y1\ =\ YCOORDS();\\
273 \ \ \ \ std::vector<{}\textbf{double}>{}::\textbf{const}\underline\ iterator\ viter;\\
274 \ \ \ \ std::ofstream\ outfile("{}matlab\underline\ out.m"{},std::ios::out);\\
275 \ \ \ \ \textbf{if}(!outfile)\ std::cerr\ <{}<{}\ "{}Unable\ to\ open\ the\ matlab\ outfile$\backslash$n"{};\\
276 \ \ \ \ outfile\ <{}<{}\ "{}clear;$\backslash$n"{};\\
277 \ \ \ \ outfile\ <{}<{}\ "{}holdon=ishold;$\backslash$n"{};\\
278 \ \ \ \ \\
279 \ \ \ \ \textbf{for}(\textbf{int}\ j\ =\ 0\ ;\ j\ <{}\ ydim\ ;\ ++j)\{\\
280 \ \ \textbf{for}(\textbf{int}\ i\ =\ 0\ ;\ i\ <{}\ xdim\ ;\ ++i)\{\\
281 \ \ \ \ \ \ \textbf{int}\ no\ =\ i\ +\ j$\ast$xdim;\\
282 \ \ \ \ \ \ outfile\ <{}<{}\ "{}x1("{}\ <{}<{}\ i+1\ <{}<{}\ "{},"{}\ <{}<{}\ j+1\ <{}<{}\ "{})="{}\ <{}<{}\ x1[no]\ <{}<{}\ "{};\ \ \ "{}\\
283 \ \ \ \ \ \ \ \ <{}<{}\ "{}y1("{}\ <{}<{}\ i+1\ <{}<{}\ "{},"{}\ <{}<{}\ j+1\ <{}<{}\ "{})="{}\ <{}<{}\ y1[no]\ <{}<{}\ "{};"{}\ <{}<{}\ std::endl;\\
284 \ \ \}\\
285 \ \ \ \ \}\\
286 \ \ \ \ \\
287 \ \ \ \ outfile\ <{}<{}\ "{}m\ =\ \ "{}\ <{}<{}\ xdim\ <{}<{}\ std::endl;\\
288 \ \ \ \ outfile\ <{}<{}\ "{}n\ =\ \ "{}\ <{}<{}\ ydim\ <{}<{}\ std::endl;\\
289 \ \ \ \ \\
290 \ \ \ \ outfile\ <{}<{}\ "{}plot(x1(1,:),y1(1,:),'r');\ hold\ on"{}\ <{}<{}\ std::endl;\\
291 \ \ \ \ outfile\ <{}<{}\ "{}plot(x1(m,:),y1(m,:),'r');"{}\ <{}<{}\ std::endl;\\
292 \ \ \ \ outfile\ <{}<{}\ "{}plot(x1(:,1),y1(:,1),'r');"{}\ <{}<{}\ std::endl;\\
293 \ \ \ \ outfile\ <{}<{}\ "{}plot(x1(:,n),y1(:,n),'r');"{}\ <{}<{}\ std::endl;\\
294 \ \ \ \ \\
295 \ \ \ \ outfile\ <{}<{}\ "{}\%\ Plot\ internal\ grid\ lines$\backslash$n"{};\\
296 \ \ \ \ outfile\ <{}<{}\ "{}\textbf{for}\ i=2:m-{}1,\ plot(x1(i,:),y1(i,:),'b');\ end$\backslash$n"{};\\
297 \ \ \ \ outfile\ <{}<{}\ "{}\textbf{for}\ j=2:n-{}1,\ plot(x1(:,j),y1(:,j),'b');\ end$\backslash$n"{};\\
298 \ \ \ \ \ \ \ \ \\
299 \ \ \ \ outfile\ <{}<{}\ "{}\textbf{if}\ (\textasciitilde holdon),\ hold\ off,\ end"{}\ <{}<{}\ std::endl;\\
300 \ \ \ \ \\
301 \ \ \ \ outfile\ <{}<{}\ "{}axis\ off;$\backslash$n"{};\\
302 \ \ \ \ \\
303 \ \ \ \ outfile.close();\\
304 \}\\
305 \ \\
306 \textbf{void}\ Domain::GMV\underline\ Writer(std::ofstream\ \&\ outFile)\{\\
307 \ \ \ \ \\
308 \ \ \ \ std::vector<{}\textbf{double}>{}\ xcoords\ =\ XCOORDS();\\
309 \ \ \ \ std::vector<{}\textbf{double}>{}\ ycoords\ =\ YCOORDS();\\
310 \ \ \ \ \ \ \ \ \\
311 \ \ \ \ outFile\ <{}<{}\ "{}gmvinput\ ascii$\backslash$n"{};\\
312 \ \ \ \ outFile\ <{}<{}\ "{}nodes\ \ "{}\ <{}<{}\ xdim$\ast$ydim\ <{}<{}\ std::endl;\\
313 \ \ \ \ \\
314 \ \ \ \ \textbf{for}(\textbf{int}\ j\ =\ 0\ ;\ j\ <{}\ ydim\ ;\ ++j)\{\\
315 \ \ \textbf{for}(\textbf{int}\ i\ =\ 0\ ;\ i\ <{}\ xdim\ ;\ ++i)\{\\
316 \ \ \ \ \ \ \textbf{int}\ no\ =\ i\ +\ j$\ast$xdim;\\
317 \ \ \ \ \ \ outFile\ <{}<{}\ xcoords[no]\ <{}<{}\ "{}\ \ \ \ \ \ \ \ \ "{};\\
318 \ \ \ \ \ \ \\
319 \ \ \}\\
320 \ \ \ \ \}\\
321 \ \ \ \ outFile\ <{}<{}\ std::endl\ <{}<{}\ std::endl;\\
322 \ \ \ \ \textsl{//writing\ y\ coord}\\
323 \ \ \ \ \textbf{for}(\textbf{int}\ j\ =\ 0\ ;\ j\ <{}\ ydim\ ;\ ++j)\{\\
324 \ \ \textbf{for}(\textbf{int}\ i\ =\ 0\ ;\ i\ <{}\ xdim\ ;\ ++i)\{\\
325 \ \ \ \ \ \ \textbf{int}\ no\ =\ i\ +\ j$\ast$xdim;\\
326 \ \ \ \ \ \ outFile\ <{}<{}\ ycoords[no]\ <{}<{}\ "{}\ \ \ \ \ \ \ \ "{}\ ;\\
327 \ \ \ \ \ \ \\
328 \ \ \}\\
329 \ \ \ \ \}\\
330 \ \ \ \ outFile\ <{}<{}\ std::endl\ <{}<{}\ std::endl;\\
331 \ \ \ \ \textsl{//forming\ cells}\\
332 \ \ \ \ outFile\ <{}<{}\ "{}cells\ \ "{}\ <{}<{}\ (xdim-{}1)$\ast$(ydim-{}1)\ <{}<{}\ std::endl;\\
333 \ \ \ \ \textbf{for}(\textbf{int}\ j\ =\ 0\ ;\ j\ <{}\ (ydim-{}1)\ ;\ ++j)\{\\
334 \ \ \textbf{for}(\textbf{int}\ i\ =\ 0\ ;\ i\ <{}\ (xdim-{}1)\ ;\ ++i)\{\\
335 \ \ \ \ \ \ \textbf{int}\ no\ =\ \ (\ i\ +\ j$\ast$(xdim)\ \ )+1;\\
336 \ \ \ \ \ \ \textbf{int}\ no1\ =\ (\ i\ +\ (j+1)$\ast$(xdim)+1);\\
337 \ \ \ \ \ \ outFile\ <{}<{}\ "{}quad\ \ 4\ \ "{}\ <{}<{}\ std::endl;\\
338 \ \ \ \ \ \ outFile\ <{}<{}\ no\ <{}<{}\ "{}\ \ \ "{}\ <{}<{}\ no+1\ <{}<{}\ "{}\ \ \ "{}\\
339 \ \ \ \ \ \ \ \ <{}<{}\ no1+1\ <{}<{}\ "{}\ \ "{}\ \ <{}<{}\ no1\ <{}<{}\ std::endl;\\
340 \ \ \}\\
341 \ \ \ \ \}\\
342 \ \ \ \ outFile\ <{}<{}\ std::endl;\\
343 \ \ \ \ \ \ \ \ \\
344 \ \ \ \ outFile\ <{}<{}\ std::endl\ <{}<{}\ "{}endgmv$\backslash$n"{};\\
345 \ \ \ \ outFile.close();\\
346 \}\\
347 \textsl{//+++++++++++++}\\
348 Matrix\ Domain::MeshX()\{\\
349 \ \ \ \ \textbf{return}\ x;\\
350 \}\\
351 Matrix\ Domain::MeshY()\{\\
352 \ \ \ \ \\
353 \ \ \ \ \textbf{return}\ y;\\
354 \}\\
355 \textsl{//++++++++++++++++++++}\\
356 \textbf{double}\ Domain::G22(\textbf{unsigned}\ \textbf{int}\ i\ ,\ \textbf{unsigned}\ \textbf{int}\ j)\{\\
357 \ \ \ \ \textbf{double}\ x1;\textbf{double}\ x2;\textbf{double}\ y1;\textbf{double}\ y2;\\
358 \ \ \ \ x1\ =\ x(i,j-{}1);\ x2\ =\ x(i,j+1)\ ;\\
359 \ \ \ \ y1\ =\ y(i,j-{}1);\ y2\ =\ y(i,j+1);\\
360 \ \ \ \ \textbf{double}\ g22\ =\ std::pow((x2-{}x1)/(2.0$\ast$del\underline\ eta),2)\ +\\
361 \ \ std::pow((y2-{}y1)/(2.0$\ast$del\underline\ eta),2);\\
362 \ \ \ \ \textbf{return}\ g22;\\
363 \}\\
364 \textbf{double}\ Domain::G11(\textbf{unsigned}\ \textbf{int}\ i\ ,\ \textbf{unsigned}\ \textbf{int}\ j)\{\\
365 \ \ \ \ \textbf{double}\ x1\ =\ x(i-{}1,j);\ \textbf{double}\ x2\ =\ x(i+1,j);\\
366 \ \ \ \ \textbf{double}\ y1\ =\ y(i-{}1,j);\ \textbf{double}\ y2\ =\ y(i+1,j);\\
367 \ \ \ \ \textbf{double}\ g11\ =\ std::pow((x2-{}x1)/(2.0$\ast$del\underline\ xi),2)\ +\\
368 \ \ std::pow((y2-{}y1)/(2.0$\ast$del\underline\ xi),2);\\
369 \ \ \ \ \textbf{return}\ g11;\\
370 \}\\
371 \textbf{double}\ Domain::X\underline\ xi(\textbf{unsigned}\ \textbf{int}\ i\ ,\ \textbf{unsigned}\ \textbf{int}\ j)\{\\
372 \ \ \ \ \textbf{double}\ x\underline\ xi;\\
373 \ \ \ \ x\underline\ xi\ =\ (x(i+1,j)-{}x(i-{}1,j))/(2.0$\ast$del\underline\ xi);\\
374 \ \ \ \ \textbf{return}\ x\underline\ xi;\\
375 \}\\
376 \textbf{double}\ Domain::X\underline\ eta(\textbf{unsigned}\ \textbf{int}\ i\ ,\ \textbf{unsigned}\ \textbf{int}\ j)\{\\
377 \ \ \ \ \textbf{double}\ x\underline\ eta;\\
378 \ \ \ \ x\underline\ eta\ =\ (x(i,j+1)-{}x(i,j-{}1))/(2.0$\ast$del\underline\ eta);\\
379 \ \ \ \ \textbf{return}\ x\underline\ eta;\\
380 \}\\
381 \textbf{double}\ Domain::Y\underline\ xi(\textbf{unsigned}\ \textbf{int}\ i\ ,\ \textbf{unsigned}\ \textbf{int}\ j)\{\\
382 \ \ \ \ \textbf{double}\ y\underline\ xi;\\
383 \ \ \ \ y\underline\ xi\ =\ (y(i+1,j)-{}y(i-{}1,j))/(2.0$\ast$del\underline\ xi);\\
384 \ \ \ \ \textbf{return}\ y\underline\ xi;\\
385 \}\\
386 \textbf{double}\ Domain::Y\underline\ eta(\textbf{unsigned}\ \textbf{int}\ i\ ,\ \textbf{unsigned}\ \textbf{int}\ j)\{\\
387 \ \ \ \ \textbf{double}\ y\underline\ eta;\\
388 \ \ \ \ y\underline\ eta\ =\ (y(i,j+1)-{}y(i,j-{}1))/(2.0$\ast$del\underline\ eta);\\
389 \ \ \ \ \textbf{return}\ y\underline\ eta;\\
390 \}\\
391 \textbf{double}\ Domain::X\underline\ xieta(\textbf{unsigned}\ \textbf{int}\ i\ ,\ \textbf{unsigned}\ \textbf{int}\ j)\{\\
392 \ \ \ \ \textbf{return}\ (x(i+1,j+1)+x(i-{}1,j-{}1)-{}x(i-{}1,j+1)-{}x(i+1,j-{}1))/(4.0$\ast$del\underline\ xi$\ast$del\underline\ eta);\\
393 \}\\
394 \textbf{double}\ Domain::Y\underline\ xieta(\textbf{unsigned}\ \textbf{int}\ i\ ,\ \textbf{unsigned}\ \textbf{int}\ j)\{\\
395 \ \ \ \ \textbf{return}\ (y(i+1,j+1)+y(i-{}1,j-{}1)-{}y(i-{}1,j+1)-{}y(i+1,j-{}1))/(4.0$\ast$del\underline\ xi$\ast$del\underline\ eta);\\
396 \}\\
397 \textbf{void}\ Domain::Fill\underline\ del\underline\ xi\underline\ eta(\textbf{double}\ xi,\textbf{double}\ eta)\{\\
398 \ \ \ \ del\underline\ xi\ =\ xi;\ del\underline\ eta\ =\ eta;\\
399 \}\\
400 std::vector<{}\textbf{double}>{}\ Domain::P11(\textbf{unsigned}\ \textbf{int}\ i\ ,\ \textbf{unsigned}\ \textbf{int}\ j)\{\\
401 \ \ \ \ \textsl{//x-{}t\ coordinate}\\
402 \ \ \ \ \textsl{//compute\ the\ jacobian\ at\ the\ point}\\
403 \ \ \ \ \textbf{double}\ s\underline\ xi\ \ \ =\ (x(i+1,j)-{}x(i-{}1,j))/(2.0$\ast$del\underline\ xi);\\
404 \ \ \ \ \textbf{double}\ s\underline\ eta\ \ =\ (x(i,j+1)-{}x(i,j-{}1))/(2.0$\ast$del\underline\ eta);\\
405 \ \ \ \ \textbf{double}\ t\underline\ xi\ \ \ =\ (y(i+1,j)-{}y(i-{}1,j))/(2.0$\ast$del\underline\ xi);\\
406 \ \ \ \ \textbf{double}\ t\underline\ eta\ \ =\ (y(i,j+1)-{}y(i,j-{}1))/(2.0$\ast$del\underline\ eta);\\
407 \ \ \ \ \textbf{double}\ det\ =\ s\underline\ xi$\ast$t\underline\ eta-{}t\underline\ xi$\ast$s\underline\ eta;\\
408 \ \ \ \ \textbf{double}\ TI\underline\ 11\ =\ \ t\underline\ eta/det;\ \textbf{double}\ \ TI\underline\ 12\ =\ -{}t\underline\ xi/det;\\
409 \ \ \ \ \textbf{double}\ TI\underline\ 21\ =\ -{}s\underline\ eta/det;\ \textbf{double}\ \ TI\underline\ 22\ =\ s\underline\ xi/det;\\
410 \ \ \ \ \textbf{double}\ s\underline\ xixi\ =\ (x(i+1,j)-{}2.0$\ast$x(i,j)+x(i-{}1,j))/(del\underline\ xi$\ast$del\underline\ xi);\\
411 \ \ \ \ \textbf{double}\ s\underline\ etaeta\ =\ (x(i,j+1)-{}2.0$\ast$x(i,j)+x(i,j-{}1))/(del\underline\ eta$\ast$del\underline\ eta);\\
412 \ \ \ \ \textbf{double}\ s\underline\ xieta\ =\ (x(i+1,j+1)+x(i-{}1,j-{}1)-{}x(i-{}1,j+1)-{}x(i+1,j-{}1))/(4.0$\ast$del\underline\ xi$\ast$del\underline\ eta);\\
413 \ \ \ \ \textbf{double}\ t\underline\ xieta\ =\ (y(i+1,j+1)+y(i-{}1,j-{}1)-{}y(i-{}1,j+1)-{}y(i+1,j-{}1))/(4.0$\ast$del\underline\ xi$\ast$del\underline\ eta);\\
414 \ \ \ \ \textbf{double}\ t\underline\ xixi\ \ \ =\ (y(i+1,j)-{}2.0$\ast$y(i,j)+y(i-{}1,j))/(del\underline\ xi$\ast$del\underline\ xi);\\
415 \ \ \ \ \textbf{double}\ t\underline\ etaeta\ =\ (y(i,j+1)-{}2.0$\ast$y(i,j)+y(i,j-{}1))/(del\underline\ eta$\ast$del\underline\ eta);\\
416 \ \ \ \ std::vector<{}\textbf{double}>{}\ P11(2);\\
417 \ \ \ \ \\
418 \ \ \ \ P11[0]\ =\ -{}(s\underline\ xixi$\ast$TI\underline\ 11+t\underline\ xixi$\ast$TI\underline\ 12);\\
419 \ \ \ \ P11[1]\ =\ -{}(s\underline\ xixi$\ast$TI\underline\ 21+t\underline\ xixi$\ast$TI\underline\ 22);\\
420 \ \ \ \ \\
421 \ \ \ \ \textbf{return}\ P11;\\
422 \}\\
423 std::vector<{}\textbf{double}>{}\ Domain::P22(\textbf{unsigned}\ \textbf{int}\ i\ ,\ \textbf{unsigned}\ \textbf{int}\ j)\{\\
424 \ \ \ \ \textsl{//x-{}t\ coordinate}\\
425 \ \ \ \ \textsl{//compute\ the\ jacobian\ at\ the\ point}\\
426 \ \ \ \ \textbf{double}\ s\underline\ xi\ \ \ =\ (x(i+1,j)-{}x(i-{}1,j))/(2.0$\ast$del\underline\ xi);\\
427 \ \ \ \ \textbf{double}\ s\underline\ eta\ \ =\ (x(i,j+1)-{}x(i,j-{}1))/(2.0$\ast$del\underline\ eta);\\
428 \ \ \ \ \textbf{double}\ t\underline\ xi\ \ \ =\ (y(i+1,j)-{}y(i-{}1,j))/(2.0$\ast$del\underline\ xi);\\
429 \ \ \ \ \textbf{double}\ t\underline\ eta\ \ =\ (y(i,j+1)-{}y(i,j-{}1))/(2.0$\ast$del\underline\ eta);\\
430 \ \ \ \ \textbf{double}\ det\ =\ s\underline\ xi$\ast$t\underline\ eta-{}t\underline\ xi$\ast$s\underline\ eta;\\
431 \ \ \ \ \textbf{double}\ TI\underline\ 11\ =\ \ t\underline\ eta/det;\ \textbf{double}\ \ TI\underline\ 12\ =\ -{}t\underline\ xi/det;\\
432 \ \ \ \ \textbf{double}\ TI\underline\ 21\ =\ -{}s\underline\ eta/det;\ \textbf{double}\ \ TI\underline\ 22\ =\ s\underline\ xi/det;\\
433 \ \ \ \ \textbf{double}\ s\underline\ xixi\ =\ (x(i+1,j)-{}2.0$\ast$x(i,j)+x(i-{}1,j))/(del\underline\ xi$\ast$del\underline\ xi);\\
434 \ \ \ \ \textbf{double}\ s\underline\ etaeta\ =\ (x(i,j+1)-{}2.0$\ast$x(i,j)+x(i,j-{}1))/(del\underline\ eta$\ast$del\underline\ eta);\\
435 \ \ \ \ \textbf{double}\ s\underline\ xieta\ =\ (x(i+1,j+1)+x(i-{}1,j-{}1)-{}x(i-{}1,j+1)-{}x(i+1,j-{}1))/(4.0$\ast$del\underline\ xi$\ast$del\underline\ eta);\\
436 \ \ \ \ \textbf{double}\ t\underline\ xieta\ =\ (y(i+1,j+1)+y(i-{}1,j-{}1)-{}y(i-{}1,j+1)-{}y(i+1,j-{}1))/(4.0$\ast$del\underline\ xi$\ast$del\underline\ eta);\\
437 \ \ \ \ \textbf{double}\ t\underline\ xixi\ \ \ =\ (y(i+1,j)-{}2.0$\ast$y(i,j)+y(i-{}1,j))/(del\underline\ xi$\ast$del\underline\ xi);\\
438 \ \ \ \ \textbf{double}\ t\underline\ etaeta\ =\ (y(i,j+1)-{}2.0$\ast$y(i,j)+y(i,j-{}1))/(del\underline\ eta$\ast$del\underline\ eta);\\
439 \ \ \ \ std::vector<{}\textbf{double}>{}\ P22(2);\\
440 \ \ \ \ P22[0]\ =\ -{}(s\underline\ etaeta$\ast$TI\underline\ 11+t\underline\ etaeta$\ast$TI\underline\ 12);\\
441 \ \ \ \ P22[1]\ =\ -{}(s\underline\ etaeta$\ast$TI\underline\ 21+t\underline\ etaeta$\ast$TI\underline\ 22);\\
442 \ \ \ \ \textbf{return}\ P22;\\
443 \}\\
444 std::vector<{}\textbf{double}>{}\ Domain::P12(\textbf{unsigned}\ \textbf{int}\ i\ ,\ \textbf{unsigned}\ \textbf{int}\ j)\{\\
445 \ \ \ \ \textsl{//x-{}t\ coordinate}\\
446 \ \ \ \ \textsl{//compute\ the\ jacobian\ at\ the\ point}\\
447 \ \ \ \ \textbf{double}\ s\underline\ xi\ \ \ =\ (x(i+1,j)-{}x(i-{}1,j))/(2.0$\ast$del\underline\ xi);\\
448 \ \ \ \ \textbf{double}\ s\underline\ eta\ \ =\ (x(i,j+1)-{}x(i,j-{}1))/(2.0$\ast$del\underline\ eta);\\
449 \ \ \ \ \textbf{double}\ t\underline\ xi\ \ \ =\ (y(i+1,j)-{}y(i-{}1,j))/(2.0$\ast$del\underline\ xi);\\
450 \ \ \ \ \textbf{double}\ t\underline\ eta\ \ =\ (y(i,j+1)-{}y(i,j-{}1))/(2.0$\ast$del\underline\ eta);\\
451 \ \ \ \ \textbf{double}\ det\ =\ s\underline\ xi$\ast$t\underline\ eta-{}t\underline\ xi$\ast$s\underline\ eta;\\
452 \ \ \ \ \textbf{double}\ TI\underline\ 11\ =\ \ t\underline\ eta/det;\ \textbf{double}\ \ TI\underline\ 12\ =\ -{}t\underline\ xi/det;\\
453 \ \ \ \ \textbf{double}\ TI\underline\ 21\ =\ -{}s\underline\ eta/det;\ \textbf{double}\ \ TI\underline\ 22\ =\ s\underline\ xi/det;\\
454 \ \ \ \ \textbf{double}\ s\underline\ xixi\ =\ (x(i+1,j)-{}2.0$\ast$x(i,j)+x(i-{}1,j))/(del\underline\ xi$\ast$del\underline\ xi);\\
455 \ \ \ \ \textbf{double}\ s\underline\ etaeta\ =\ (x(i,j+1)-{}2.0$\ast$x(i,j)+x(i,j-{}1))/(del\underline\ eta$\ast$del\underline\ eta);\\
456 \ \ \ \ \textbf{double}\ s\underline\ xieta\ =\ (x(i+1,j+1)+x(i-{}1,j-{}1)-{}x(i-{}1,j+1)-{}x(i+1,j-{}1))/(4.0$\ast$del\underline\ xi$\ast$del\underline\ eta);\\
457 \ \ \ \ \textbf{double}\ t\underline\ xieta\ =\ (y(i+1,j+1)+y(i-{}1,j-{}1)-{}y(i-{}1,j+1)-{}y(i+1,j-{}1))/(4.0$\ast$del\underline\ xi$\ast$del\underline\ eta);\\
458 \ \ \ \ \textbf{double}\ t\underline\ xixi\ \ \ =\ (y(i+1,j)-{}2.0$\ast$y(i,j)+y(i-{}1,j))/(del\underline\ xi$\ast$del\underline\ xi);\\
459 \ \ \ \ \textbf{double}\ t\underline\ etaeta\ =\ (y(i,j+1)-{}2.0$\ast$y(i,j)+y(i,j-{}1))/(del\underline\ eta$\ast$del\underline\ eta);\\
460 \ \ \ \ std::vector<{}\textbf{double}>{}\ P12(2);\\
461 \ \ \ \ P12[0]\ =\ -{}(s\underline\ xieta$\ast$TI\underline\ 11+t\underline\ xieta$\ast$TI\underline\ 12);\\
462 \ \ \ \ P12[1]\ =\ -{}(s\underline\ xieta$\ast$TI\underline\ 21+t\underline\ xieta$\ast$TI\underline\ 22);\\
463 \ \ \ \ \textbf{return}\ P12;\\
464 \}\\
465 \ \\
466  }
\normalfont\normalsize

\subsection{matrix.cpp}
{\ttfamily \raggedright \tiny
001 \#ifndef MATRIX\underline\ H\\
002 \#define MATRIX\underline\ H\\
003 \textbf{class}\ Matrix\{\\
004 \textbf{public}:\\
005 \ \ \ \ \textbf{unsigned}\ \textbf{int}\ nx,ny;\\
006 \ \ \ \ std::vector<{}std::vector<{}\textbf{double}>{}\ >{}\ Elements;\\
007 \ \ \ \ Matrix()\{\\
008 \ \ \ \ \}\\
009 \ \ \ \ \\
010 \ \ \ \ Matrix(\textbf{unsigned}\ \textbf{int}\ nx1,\textbf{unsigned}\ \textbf{int}\ ny1)\{\\
011 \ \ nx\ =\ nx1\ ;\ ny\ =\ ny1;\\
012 \ \ \textsl{//number\ of\ rows}\\
013 \ \ Elements.resize(nx);\\
014 \ \ \textsl{//fill\ each\ rows}\\
015 \ \ \textbf{for}(\textbf{unsigned}\ \textbf{int}\ i\ =\ 0\ ;\ i\ <{}\ nx\ ;\ ++i)\\
016 \ \ \ \ \ \ Elements[i].resize(ny,99.0);\\
017 \ \ \ \ \}\\
018 \ \ \ \ \\
019 \ \ \ \ \textbf{void}\ Clear()\{\\
020 \ \ \textbf{for}(\textbf{unsigned}\ \textbf{int}\ i\ =\ 0\ ;\ i\ <{}\ nx\ ;\ ++i)\\
021 \ \ \ \ \ \ Elements[i].clear();\\
022 \ \ \ \ \}\\
023 \ \\
024 \ \ \ \ \textbf{double}\&\ \textbf{operator}()(\textbf{unsigned}\ \textbf{int}\ ix,\ \textbf{unsigned}\ \textbf{int}\ iy)\{\\
025 \ \ \textbf{return}\ Elements[ix][iy];\\
026 \ \ \ \ \}\\
027 \};\\
028 \#endif\\
029 \ \\
030  }
\normalfont\normalsize

\subsection{sor\_solver.cpp}
\label{subsec_SOR}
{\ttfamily \raggedright \tiny
001 \#ifndef SOR\underline\ SOLVER\\
002 \#define SOR\underline\ SOLVER\\
003 \ \\
004 \#include\ <{}vector>{}\\
005 \ \\
006 \#include\ "{}domain.h"{}\\
007 \#include\ "{}matrix.h"{}\\
008 \ \\
009 \textbf{double}\ Mesh\underline\ Residual(Matrix\ x,\ \ \ Matrix\ y\ ,\\
010 \ \ \ \ \ \ \ \ \ Matrix\ x\underline\ old,\ Matrix\ y\underline\ old,\\
011 \ \ \ \ \ \ \ \ \ \textbf{unsigned}\ \textbf{int}\ xdim,\ \textbf{unsigned}\ \textbf{int}\ ydim)\{\\
012 \ \ \ \ \textbf{double}\ resid\ =\ 0\ ;\\
013 \ \ \ \ \textbf{for}(\textbf{int}\ j\ =\ 0\ ;\ j\ <{}\ ydim\ ;\ ++j)\{\\
014 \ \ \textbf{for}(\textbf{int}\ i\ =\ 0\ ;\ i\ <{}\ xdim\ ;\ ++i)\{\\
015 \ \ \ \ \ \ \textbf{double}\ x\underline\ resd\ =\ (x(i,j)-{}x\underline\ old(i,j));\\
016 \ \ \ \ \ \ \textbf{double}\ y\underline\ resd\ =\ (y(i,j)-{}y\underline\ old(i,j));\\
017 \ \ \ \ \ \ resid\ +=\ (x\underline\ resd$\ast$x\underline\ resd+y\underline\ resd$\ast$y\underline\ resd);\\
018 \ \ \}\\
019 \ \ \ \ \}\\
020 \ \ \ \ \textbf{return}\ std::sqrt(resid);\\
021 \}\\
022 \ \\
023 \textbf{bool}\ SORSOLVER(\ Domain\&\ physical,\ Domain\ \&\ parm\ ,\ \ \ \ \ \ \ \textbf{unsigned}\ \textbf{int}\ xdim1,\textbf{unsigned}\ \textbf{int}\ ydim1)\{\\
024 \ \ \ \ \textbf{bool}\ grid\underline\ dist\ =\ \textbf{true};\\
025 \ \ \ \ \textbf{bool}\ run\underline\ ellip\ =\ \textbf{true};\\
026 \ \\
027 \ \ \ \ \textbf{double}\ tolerance\ \ =\ 1.0e-{}4;\ \ \\
028 \ \ \ \ \textbf{unsigned}\ max\underline\ iter\ =\ 100;\\
029 \ \ \ \ \textbf{double}\ w\ =\ 1.90;\\
030 \ \ \ \ \textbf{double}\ residual\ =\ 10.0;\\
031 \ \ \ \ \textbf{unsigned}\ \textbf{int}\ iter\ =\ 0;\\
032 \ \ \ \ \\
033 \ \ \ \ \textbf{unsigned}\ \textbf{int}\ xdim\ =xdim1\ ;\ \textbf{unsigned}\ \textbf{int}\ ydim\ =\ ydim1;\\
034 \ \ \ \ Matrix\ x\underline\ old(xdim,ydim),\ y\underline\ old(xdim,ydim);\\
035 \ \ \ \ \ \ \ \ \\
036 \ \ \ \ \textbf{double}\ del\underline\ xi\ =\ 1.0/double(xdim-{}1.0);\\
037 \ \ \ \ \textbf{double}\ del\underline\ eta\ =\ 1.0/double(ydim-{}1.0);\\
038 \ \ \ \ \ \ \ \ \\
039 \ \ \ \ std::ofstream\ fout("{}gauss.dat"{},std::ios::out);\\
040 \ \ \ \ \textbf{if}(!(fout.is\underline\ open()))\ std::cerr\ <{}<{}\ "{}ERROR\ :\ UNABLE\ TO\ OPEN\ THE\ FILE\ $\backslash$''gauss.dat$\backslash$''$\backslash$n"{};\\
041 \ \ \ \ \\
042 \ \ \ \ \ \textbf{if}(run\underline\ ellip)\{\\
043 \ \ \ \textbf{while}(iter\ <{}\ max\underline\ iter\ \&\ residual\ >{}\ tolerance)\{\\
044 \ \ \ \ \ \ \ iter++;\\
045 \ \ \ \ \ \ \ x\underline\ old\ =\ physical.MeshX()\ ;\ y\underline\ old\ =\ physical.MeshY();\\
046 \ \ \ \ \ \ \ \textbf{for}(\textbf{unsigned}\ \textbf{int}\ j\ =\ 1\ ;\ j\ <{}\ ydim-{}1\ ;\ ++j)\{\\
047 \ \ \ \ \ \textbf{for}(\textbf{unsigned}\ \textbf{int}\ i\ =\ 1\ ;\ i\ <{}\ xdim-{}1\ ;\ ++i)\{\\
048 \ \ \ \ \ \ \ \ \ \\
049 \ \ \ \ \ \ \ \ \textbf{double}\ g22\ \ \ \ \ =\ physical.G22(i,j);\\
050 \ \ \ \ \ \ \ \ \textbf{double}\ g11\ \ \ \ \ =\ physical.G11(i,j);\\
051 \ \ \ \ \ \ \ \ \textbf{double}\ x\underline\ xi\ \ \ \ =\ physical.X\underline\ xi(i,j);\\
052 \ \ \ \ \ \ \ \ \textbf{double}\ x\underline\ eta\ \ \ =\ physical.X\underline\ eta(i,j);\\
053 \ \ \ \ \ \ \ \ \textbf{double}\ y\underline\ xi\ \ \ \ =\ physical.Y\underline\ xi(i,j);\\
054 \ \ \ \ \ \ \ \ \textbf{double}\ y\underline\ eta\ \ \ =\ physical.Y\underline\ eta(i,j);\\
055 \ \ \ \ \ \ \ \ \textbf{double}\ x\underline\ xieta\ =\ physical.X\underline\ xieta(i,j);\\
056 \ \ \ \ \ \ \ \ \textbf{double}\ y\underline\ xieta\ =\ physical.Y\underline\ xieta(i,j);\\
057 \ \ \ \ \ \ \ \ \textbf{double}\ g12\ =\ x\underline\ xi$\ast$x\underline\ eta+y\underline\ xi$\ast$y\underline\ eta;\\
058 \ \ \ \ \ \ \ \ \textbf{double}\ g\ =\ std::pow(x\underline\ xi$\ast$y\underline\ eta-{}y\underline\ xi$\ast$x\underline\ eta,2);\\
059 \ \ \ \ \ \ \ \ \\
060 \ \ \ \ \ \ \ \ std::vector<{}\textbf{double}>{}\ P11\ ,\ P22,P12;\\
061 \ \ \ \ \ \ \ \ \\
062 \ \ \ \ \ \ \ \ \textbf{if}(grid\underline\ dist)\{\\
063 \ \ \ \ \ \ P11\ =\ parm.P11(i,j);\\
064 \ \ \ \ \ \ P22\ =\ parm.P22(i,j);\\
065 \ \ \ \ \ \ P12\ =\ parm.P12(i,j);\\
066 \ \ \ \ \ \ \ \ \}\textbf{else}\{\\
067 \ \ \ \ \ \ P11.push\underline\ back(0);P11.push\underline\ back(0);\\
068 \ \ \ \ \ \ P22.push\underline\ back(0);P22.push\underline\ back(0);\\
069 \ \ \ \ \ \ P12.push\underline\ back(0);P12.push\underline\ back(0);\\
070 \ \ \ \ \ \ \ \ \}\\
071 \ \ \ \ \ \ \ \ \\
072 \ \ \ \ \ \ \ \ \textbf{double}\ tmpx,\ tmpy;\\
073 \ \ \ \ \ \ \ \ \\
074 \ \ \ \ \ \ \ \ tmpx\ =\ (g22$\ast$P11[0]-{}2.0$\ast$g12$\ast$P12[0]+g11$\ast$P22[0])$\ast$x\underline\ xi+\\
075 \ \ \ \ \ \ (g22$\ast$P11[1]-{}2.0$\ast$g12$\ast$P12[1]+g11$\ast$P22[1])$\ast$x\underline\ eta;\\
076 \ \ \ \ \ \ \ \ \\
077 \ \ \ \ \ \ \ \ tmpy\ =\ (g22$\ast$P11[0]-{}2.0$\ast$g12$\ast$P12[0]+g11$\ast$P22[0])$\ast$y\underline\ xi+\\
078 \ \ \ \ \ \ (g22$\ast$P11[1]-{}2.0$\ast$g12$\ast$P12[1]+g11$\ast$P22[1])$\ast$y\underline\ eta;\\
079 \ \ \ \ \ \ \ \ \\
080 \ \ \ \ \ \ \ \ \textbf{double}\ lhsx\ =\ 2.0$\ast$(g22/(del\underline\ xi$\ast$del\underline\ xi)+g11/(del\underline\ eta$\ast$del\underline\ eta));\\
081 \ \ \ \ \ \ \ \ \textbf{double}\ rhsx\ =\ g22$\ast$(physical.XCOORD(i+1,j)+physical.XCOORD(i-{}1,j))/(del\underline\ xi$\ast$del\underline\ xi)\\
082 \ \ \ \ \ \ +g11$\ast$(physical.XCOORD(i,j+1)+physical.XCOORD(i,j-{}1))/(del\underline\ eta$\ast$del\underline\ eta)\\
083 \ \ \ \ \ \ -{}2.0$\ast$g12$\ast$x\underline\ xieta\ +\ tmpx;\\
084 \ \ \ \ \ \ \ \ \\
085 \ \ \ \ \ \ \ \ physical.XCOORD(i,j)\ =\ physical.XCOORD(i,j)\ +\ w$\ast$(rhsx/lhsx-{}physical.XCOORD(i,j));\ \\
086 \ \ \ \ \ \ \ \ \\
087 \ \ \ \ \ \ \ \ \textbf{double}\ lhsy\ =\ lhsx;\\
088 \ \ \ \ \ \ \ \ \textbf{double}\ rhsy\ =\ g22$\ast$(physical.YCOORD(i+1,j)+physical.YCOORD(i-{}1,j))/(del\underline\ xi$\ast$del\underline\ xi)\\
089 \ \ \ \ \ \ +g11$\ast$(physical.YCOORD(i,j+1)+physical.YCOORD(i,j-{}1))/(del\underline\ eta$\ast$del\underline\ eta)\\
090 \ \ \ \ \ \ -{}2.0$\ast$g12$\ast$y\underline\ xieta+\ tmpy;\\
091 \ \ \ \ \ \ \ \ \\
092 \ \ \ \ \ \ \ \ physical.YCOORD(i,j)\ =\ physical.YCOORD(i,j)\ +\ w$\ast$(rhsy/lhsy-{}physical.YCOORD(i,j));\\
093 \ \ \ \ \}\\
094 \ \ \ \ \ \ \}\\
095 \ \ \ \ \ \ \\
096 \ \ \ \ \ \ \textbf{double}\ mesh\underline\ resid\ =\ Mesh\underline\ Residual(physical.MeshX(),physical.MeshY(),x\underline\ old,y\underline\ old,xdim,ydim)\\
097 \ \ \ \ /((xdim-{}2)$\ast$(ydim-{}2));\\
098 \ \ \ \ \ \ std::cout\ <{}<{}\ "{}Iteration\ =\ \ \ "{}\ <{}<{}\ iter\ <{}<{}\ "{},\ \ Residual\ =\ \ "{}\\
099 \ \ \ \ \ \ \ \ \ \ <{}<{}\ mesh\underline\ resid\ <{}<{}\ std::endl;\\
100 \ \ \ \ \ \ \\
101 \ \ \ \ \ \ fout\ <{}<{}\ iter\ <{}<{}\ "{}\ \ \ \ \ \ \ \ \ "{}\ <{}<{}\ mesh\underline\ resid\ <{}<{}\ std::endl;\\
102 \ \ \ \ \ \ residual\ =\ mesh\underline\ resid;\\
103 \ \ \ \}\\
104 \ \ \ \ \ \}\\
105 \ \ \ \ \ \\
106 \ \ \ \ \ \textbf{return}\ \textbf{true};\\
107 \}\ \ \ \ \ \ \\
108 \#endif\\
109 \ \\
110  }
\normalfont\normalsize

\subsection{functions.h}
\label{subsec_fun}
{\ttfamily \raggedright \tiny
001 \#ifndef \underline\ FUNCTIONS\underline\ \\
002 \#define \underline\ FUNCTIONS\underline\ \\
003 \ \\
004 \#include\ <{}vector>{}\\
005 \#include\ <{}cmath>{}\\
006 \#include<{}iomanip>{}\\
007 \#include<{}iostream>{}\\
008 \ \\
009 \textbf{double}\ XYcircle(\textbf{double}\ x,\ \textbf{double}\ y,\textbf{unsigned}\ \textbf{int}\ x\underline\ \textbf{or}\underline\ y)\{\\
010 \ \ \ \ \textbf{double}\ r\ =\ 1.0;\\
011 \ \ \ \ \textbf{double}\ theta\ =\ 0.0;\\
012 \ \ \ \ \textbf{double}\ Pi\ =\ 4.0$\ast$atan(1.0);\\
013 \ \ \ \ \textbf{if}(0==y)\{\\
014 \ \ theta\ =\ Pi/2.0$\ast$x;\\
015 \ \ \textbf{if}(x\underline\ \textbf{or}\underline\ y\ ==\ 1)\\
016 \ \ \ \ \ \ \textbf{return}\ r$\ast$cos(theta);\\
017 \ \ \textbf{else}\\
017 \ \ \ \ \ \ \textbf{return}\ r$\ast$sin(theta);\\
018 \ \ \ \ \}\\
019 \ \ \ \ \textbf{if}(1==x)\{\\
020 \ \ theta\ =\ Pi/2+Pi/2$\ast$y;\\
021 \ \ \textbf{if}(x\underline\ \textbf{or}\underline\ y\ ==\ 1)\\
022 \ \ \ \ \ \ \textbf{return}\ r$\ast$cos(theta);\\
023 \ \ \textbf{else}\\
023 \ \ \ \ \ \ \textbf{return}\ r$\ast$sin(theta);\\
024 \ \ \ \ \}\\
025 \ \ \ \ \textbf{if}(1==y)\{\\
026 \ \ theta=Pi+Pi/2$\ast$(1.0-{}x);\\
027 \ \ \textbf{if}(x\underline\ \textbf{or}\underline\ y\ ==\ 1)\\
028 \ \ \ \ \ \ \textbf{return}\ r$\ast$cos(theta);\\
029 \ \ \textbf{else}\\
029 \ \ \ \ \ \ \textbf{return}\ r$\ast$sin(theta);\\
030 \ \ \ \ \}\\
031 \ \ \ \ \textbf{if}(0==x)\{\\
032 \ \ theta\ =\ 3.0$\ast$Pi/2.0+Pi/2.0$\ast$(1.0-{}y);\\
033 \ \ \textbf{if}(x\underline\ \textbf{or}\underline\ y\ ==\ 1)\\
034 \ \ \ \ \ \ \textbf{return}\ r$\ast$cos(theta);\\
035 \ \ \textbf{else}\\
035 \ \ \ \ \ \ \textbf{return}\ r$\ast$sin(theta);\\
036 \ \ \ \ \}\\
037 \}\\
038 \#endif\\
039 \ \\
040  }
\normalfont\normalsize

\subsection{makefile}
{\ttfamily \raggedright \tiny
001 CXX\ =\ g++\\
002 SRC\ =\ main.cpp$\backslash$\\
003 domain.cpp\\
004 OBJ\ =\ main.o$\backslash$\\
005 domain.o\\
006 ell\underline\ mesh:\ \$(OBJ)\\
007 \ \ \$(CXX)\ -{}Wall\ -{}O3\ -{}o\ ellmesh\ \$(OBJ)\ \$(LIB)\\
008 clean:\\
009 \ \ rm\ -{}f\ ellmesh\\
010 \ \ rm\ -{}f\ $\ast$.o\\
011 \ \\
012  }
\normalfont\normalsize

\section{Numerical Examples}
\subsection{Example 1}
\begin{figure}
 \begin{minipage}[t]{8cm}
 \begin{center}
 \includegraphics[width=8cm,clip]{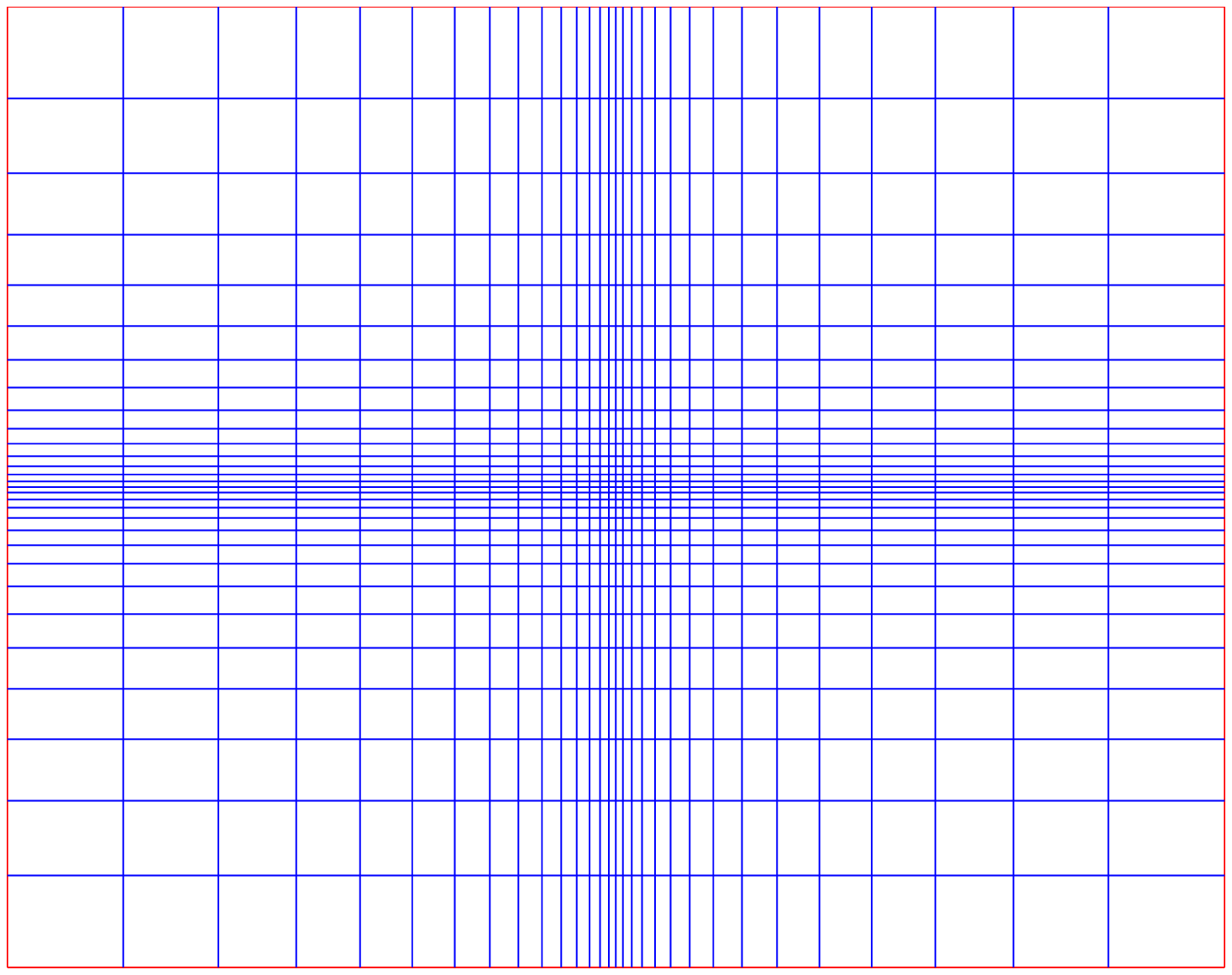}
 \caption[Transfinite Interpolation]{\label{labelFig3} Grid in the parameter space.}
 \end{center}
 \end{minipage}
 \hfill
 \begin{minipage}[t]{8cm}
 \begin{center}
 \includegraphics[width=8.0cm,clip]{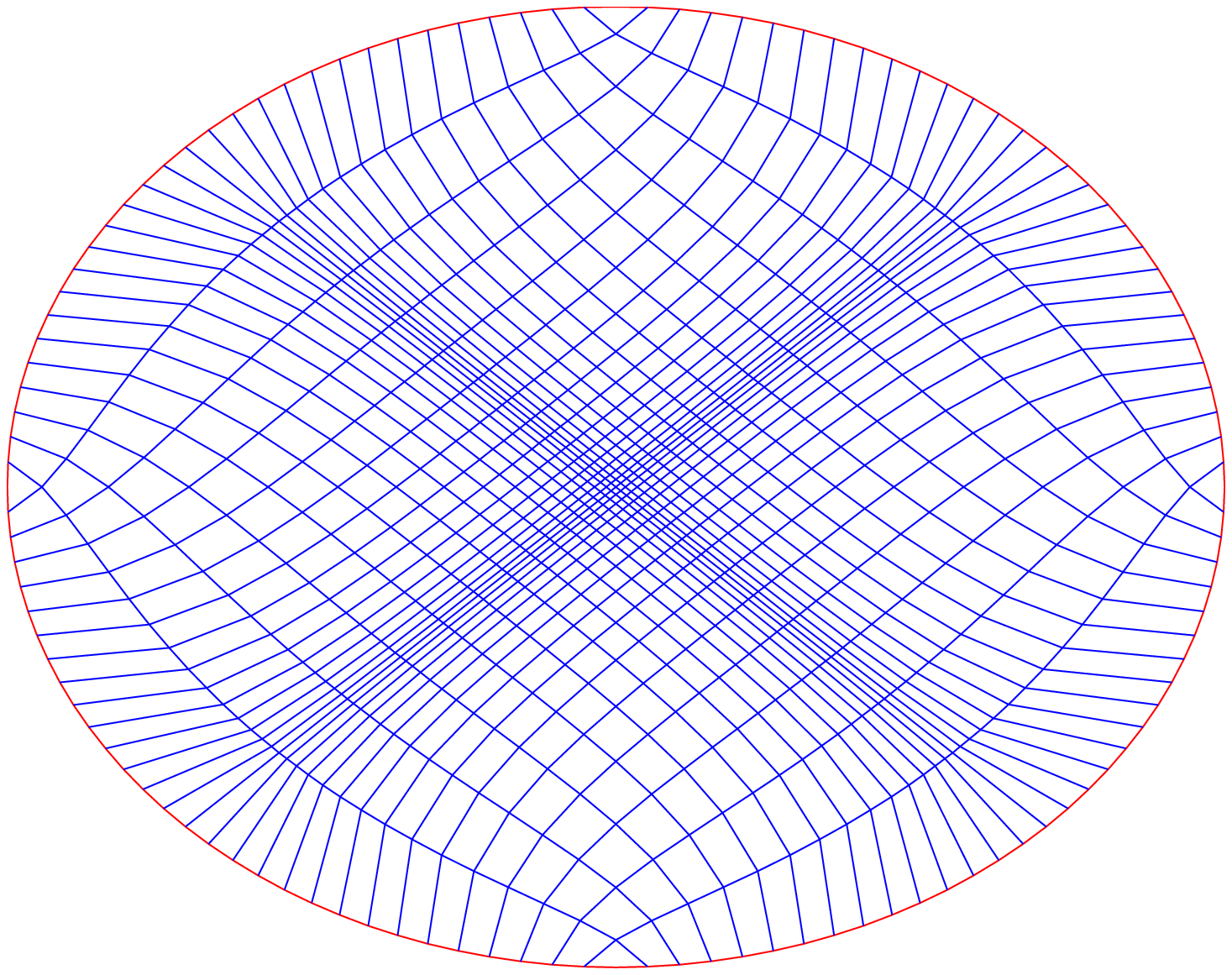}
 \caption[Elliptic Grid]{\label{labelFig4} Adapted grid by elliptic system.}
 \end{center}
 \end{minipage}
 \end{figure}

In this example, we cluster the grids along the centre lines of a circular physical domain. Figure \ref{labelFig3} shows the grid in the parameter space for concentrating grids at the centre lines of the physical space. Grid density in the physical space is determined by the grid density in the parameter space. Figure \ref{labelFig4} shows the converged grid in the physical space. For generating this grid the lines 027 and 028 in the subsection \ref{subsec_main.cpp} are used.
\subsection{Example 2}
\begin{figure}
 \begin{minipage}[t]{8cm}
 \begin{center}
 \includegraphics[width=8cm,clip]{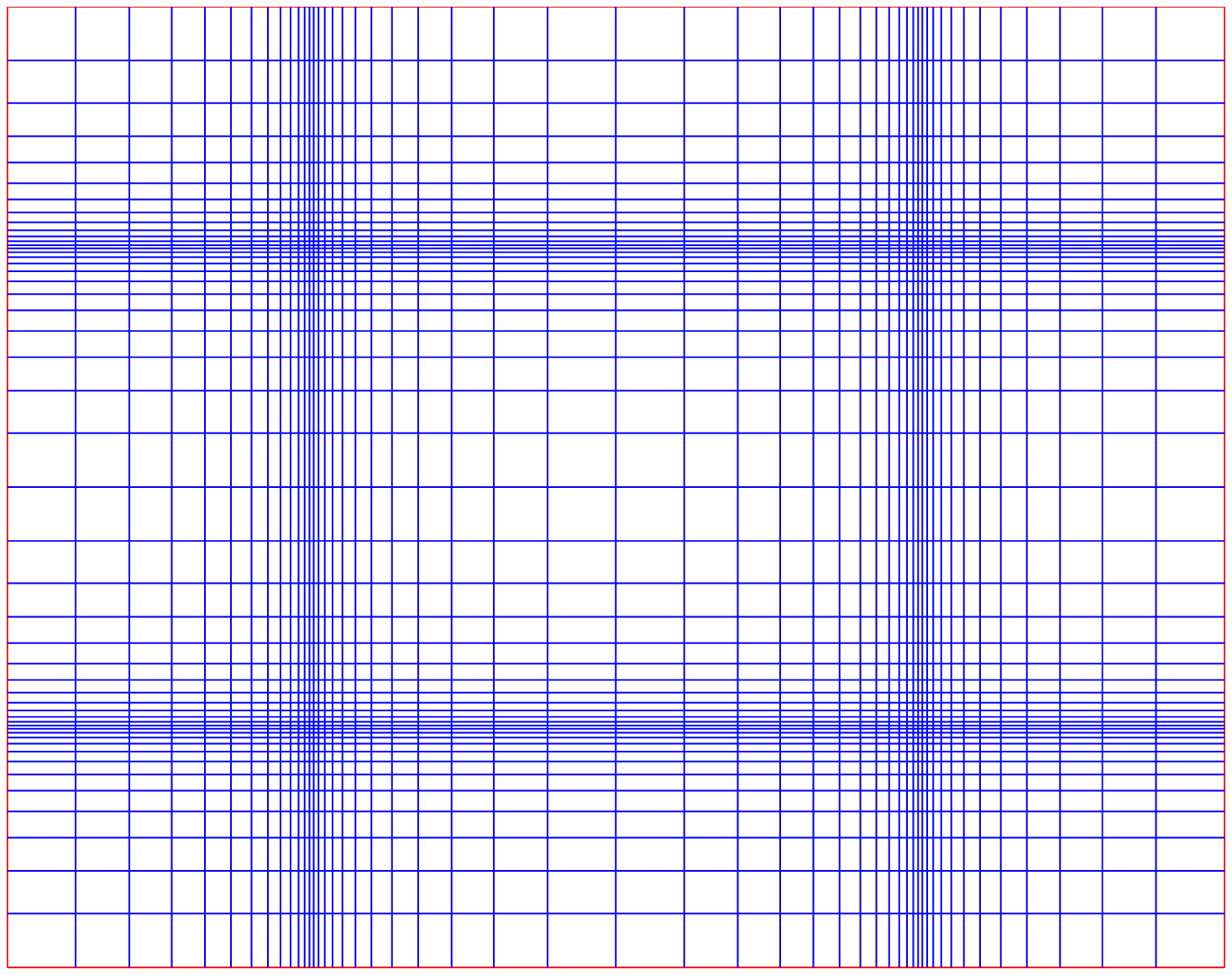}
 \caption[Transfinite Interpolation]{\label{fig:ex2_para} Grid in the parameter space.}
 \end{center}
 \end{minipage}
 \hfill
 \begin{minipage}[t]{8cm}
 \begin{center}
 \includegraphics[width=8.0cm,clip]{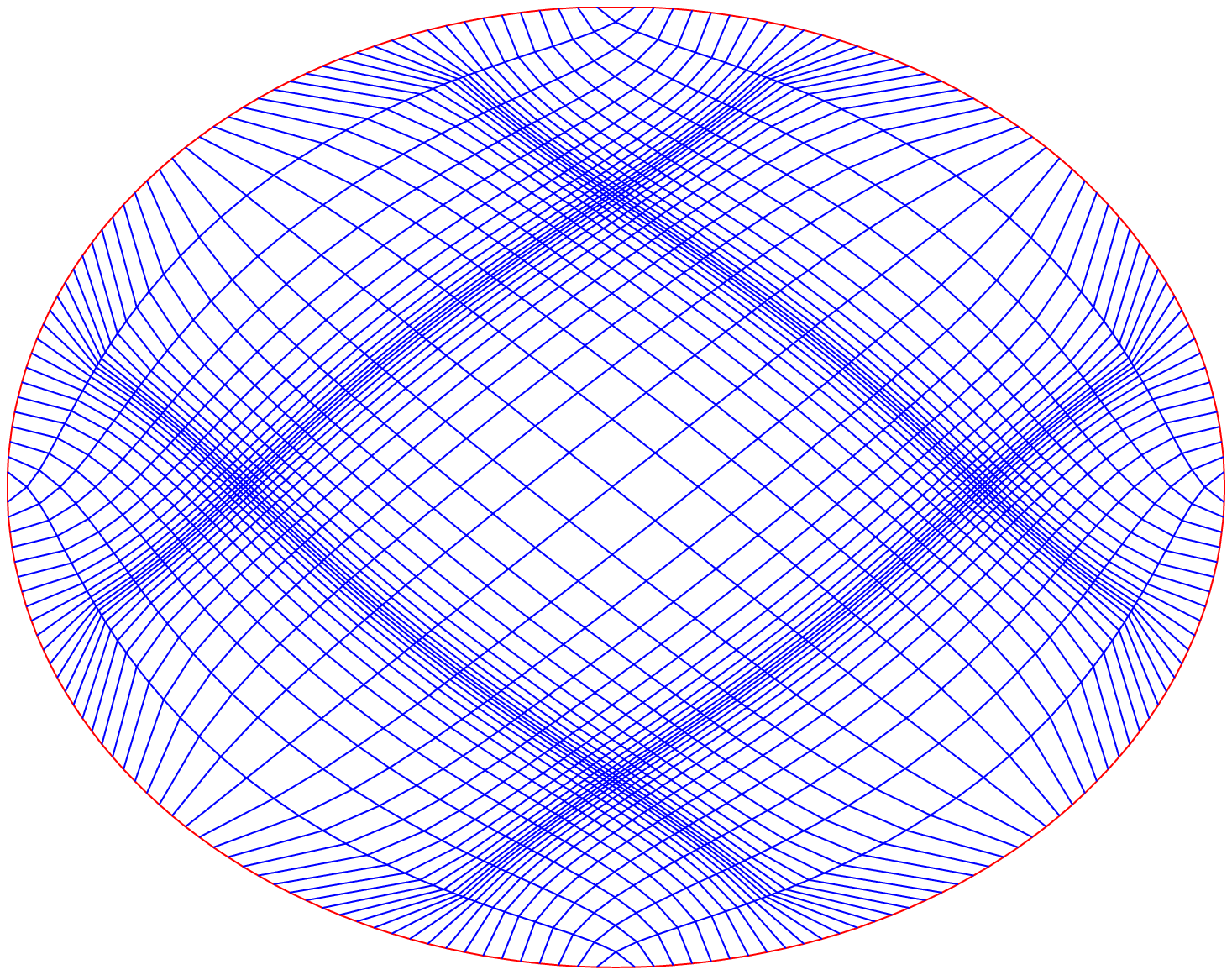}
 \caption[Elliptic Grid]{\label{fig:ex2_phy} Adapted grid by elliptic system.}
 \end{center}
 \end{minipage}
 \end{figure}

See the Figure \ref{fig:ex2_para} for grids in the parameter space and the Figure \ref{fig:ex2_phy} for the converged grids in the physical space. For generating the grids the lines 030 and 031 of the subsection \ref{subsec_main.cpp} are used.
\subsection{Example 3}
\begin{figure}
 \begin{minipage}[t]{8cm}
 \begin{center}
 \includegraphics[width=8cm,clip]{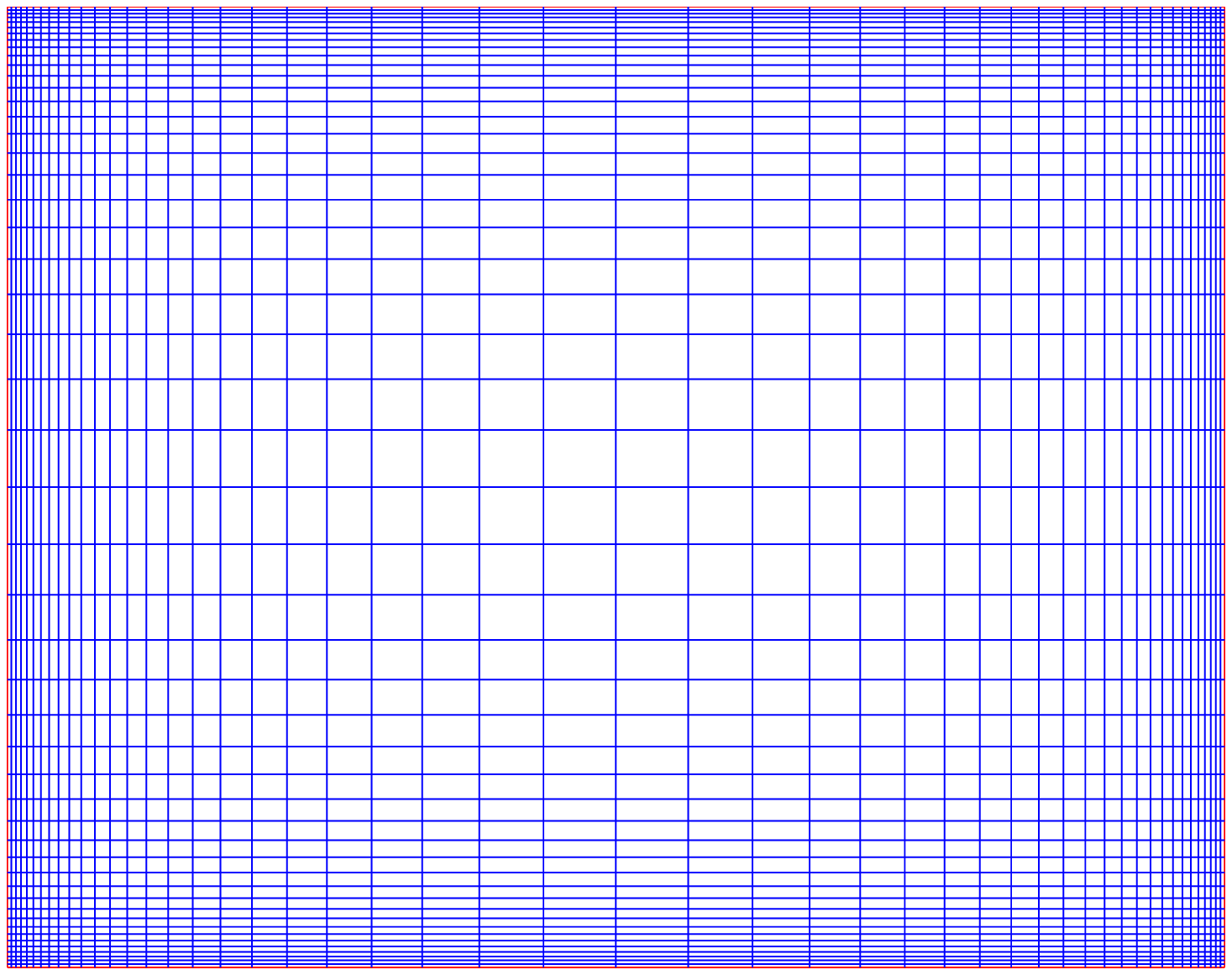}
 \caption[Transfinite Interpolation]{\label{fig:ex3_para} Grid in the parameter space.}
 \end{center}
 \end{minipage}
 \hfill
 \begin{minipage}[t]{8cm}
 \begin{center}
 \includegraphics[width=8.0cm,clip]{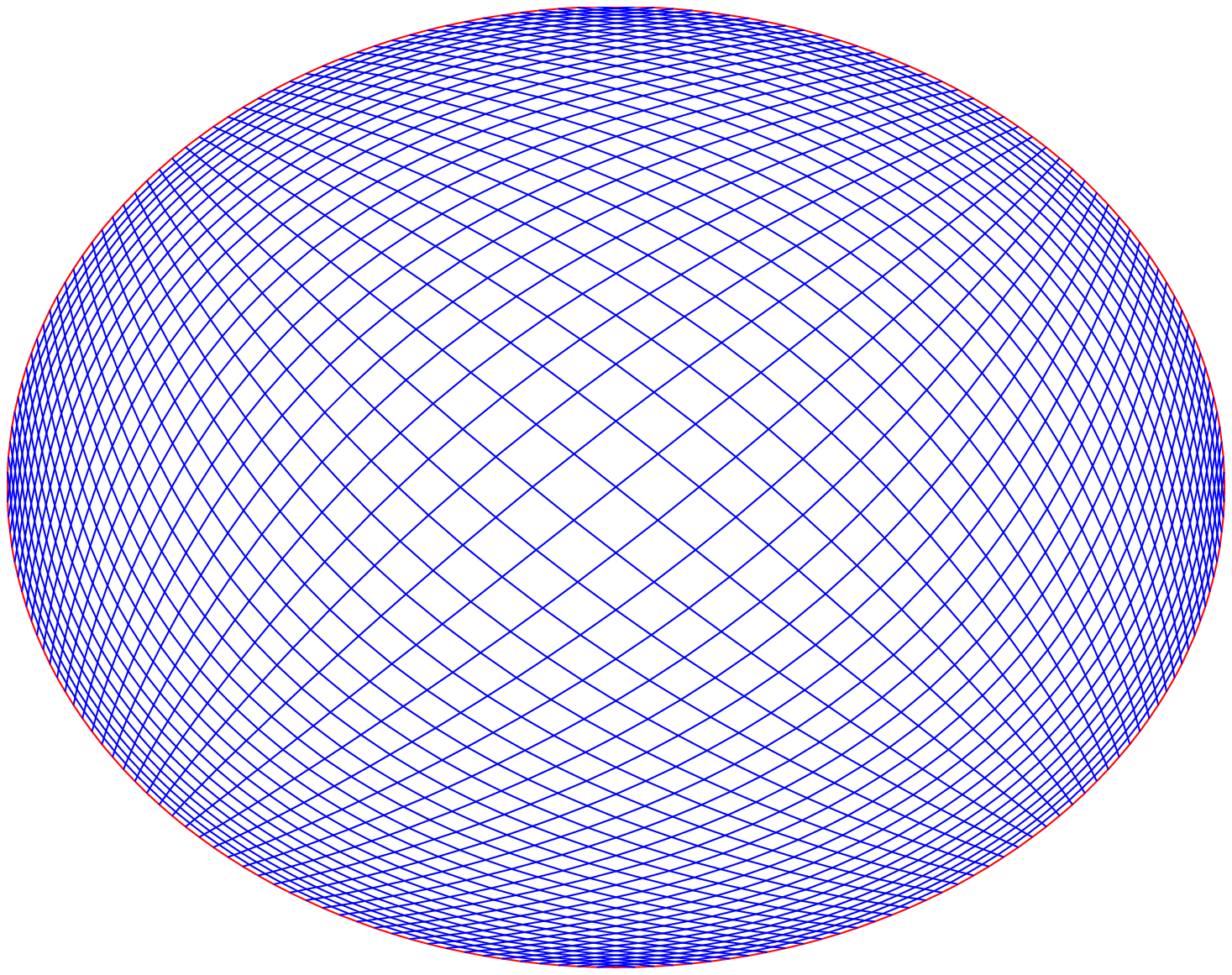}
 \caption[Elliptic Grid]{\label{fig:ex3_phy} Adapted grid by elliptic system.}
 \end{center}
 \end{minipage}
 \end{figure}

In this example, we are interested in concentrating grids at the boundary of the physical space. Figure \ref{fig:ex3_para} shows the grid in the parameter space for concentrating grids at the boundary of the physical domain. The converged grids in the physical space is shown in the Figure \ref{fig:ex3_phy}. For generating the grids the lines 033 and 034 of the subsection \ref{subsec_main.cpp} are used.
\section{Conclusions}
An elliptic system for generating adaptive quadrilateral meshes in curved domains has been presented. A C$^{++}$ implementation of the presented technique is also given. Three examples are reported for demonstrating the effectiveness of the technique and the implementation. Since, the quadrilateral meshes are very extensively used for numerical simulations and grid adaptation is required for capturing many important phenomenon such as the boundary layers. Thus, this software package can be a very useful tool. 

Our package is freely available at {\texttt{{www.mi.uib.no/$\sim$sanjay}}}. The authors want to mention that they donot know of any freely available software package that can generate adaptive quadrilateral meshes.
%
%

\end{spacing}
\end{document}